\documentclass[a4paper,11pt]{article}
\usepackage[T1]{fontenc}
\usepackage{lmodern}
\usepackage{microtype}
\usepackage[margin=24mm,top=24mm,bottom=25mm]{geometry}
\usepackage{amsmath,amssymb,amsthm,mathtools}
\usepackage{booktabs}
\usepackage{needspace}
\usepackage{xcolor}
\usepackage[unicode,colorlinks=true,linkcolor=blue!45!black,citecolor=blue!45!black,urlcolor=blue!45!black]{hyperref}
\usepackage{bookmark}
\hypersetup{pdftitle={Cohomology Operations and the Inoue-Kabaya Map via Loop Space Models},pdfauthor={Takefumi Nosaka}}
\hypersetup{pdfkeywords={quandle cohomology, Inoue–Kabaya map, rack space, James construction, half-shuffle, Pontryagin–Thomas operations}}
\newtheoremstyle{articleplain}{7pt}{7pt}{\itshape}{}{\bfseries}{.}{.6em}{}
\newtheoremstyle{articledefinition}{7pt}{7pt}{\normalfont}{}{\bfseries}{.}{.6em}{}
\theoremstyle{articleplain}
\newtheorem{theorem}{Theorem}[section]
\newtheorem{lemma}[theorem]{Lemma}
\newtheorem{proposition}[theorem]{Proposition}
\newtheorem{corollary}[theorem]{Corollary}
\theoremstyle{articledefinition}
\newtheorem{remark}[theorem]{Remark}

\newtheorem{example}[theorem]{Example}

\numberwithin{equation}{section}
\DeclareMathOperator{\As}{As}
\DeclareMathOperator{\Inn}{Inn}
\DeclareMathOperator{\Stab}{Stab}
\DeclareMathOperator{\hocolim}{hocolim}
\DeclareMathOperator{\hofib}{hofib}
\DeclareMathOperator{\Cofib}{Cofib}
\DeclareMathOperator{\im}{im}
\DeclareMathOperator{\Tor}{Tor}
\DeclareMathOperator{\Ext}{Ext}
\DeclareMathOperator{\Hom}{Hom}
\DeclareMathOperator{\id}{id}
\newcommand{\Z}{\mathbb Z}
\newcommand{\F}{\mathbb F}

\newcommand{\eps}{\varepsilon}
\newcommand{\susp}{\mathrm{s}}
\title{\bfseries  Loop Space Models for Quandle Cohomology \\
Operations and Inoue--Kabaya Map}
\author{Takefumi Nosaka}
\date{}
\begin{document}
\maketitle
\begin{abstract}
%\begin{abstract}
For connected quandles, a map from the suspended shadow rack space to the
associated orbit space realizes the reduced shadow Inoue--Kabaya map up to
integral chain homotopy. For a finite augmentation group with point stabilizer
$H$, the rack half-product agrees with the James half-shuffle over fields of
characteristic $p$ with $p\nmid |H|$. In a family of finite symplectic
quandles, we compute the orders of cup powers and Pontryagin--Thomas values
for degree-four classes in the image of the ordinary map. For each prime $p$,
these classes have unbounded nilpotence indices over $\mathbb{Z}_{(p)}$ even
at fixed quandle type.
\end{abstract}
\begin{quote}
\small
\noindent\textbf{2020 Mathematics Subject Classification.}\enspace Primary 57K12. Secondary 55P35, 55S05.\par
\noindent\textbf{Keywords.}\enspace quandle cohomology, Inoue--Kabaya map, rack space, James construction, half-shuffle, Pontryagin--Thomas operations.
\end{quote}

\section{Introduction}

We use loop spaces to study cohomology operations on the image of the Inoue--Kabaya map.
Does quandle type bound the nilpotence indices of degree-four classes in this image over $\Z_{(p)}$?

Let $H$ be a subgroup of a group $G$, and let $h\in Z(H)$. For each $y=Hk$ in the homogeneous
space $X\coloneqq H\backslash G$, right multiplication by $k^{-1}hk$ is independent of the
choice of representative $k$ and fixes $y$. Writing $x*y$ for the image of $x$ under this right
multiplication defines a quandle operation on $X$. For example, if $H$ is the centralizer of $h$
in $G$, the map $Hk\mapsto k^{-1}hk$ identifies $X$ with the conjugacy class of $h$, with
operation $x*y=y^{-1}xy$.
Every connected quandle admits a presentation $H\backslash G$ of the above form with $h\in Z(H)$. Take $G$ to be the group generated by the
right translations $S_y\colon x\mapsto x*y$, which acts transitively on $X$, and take $H$ to
be the stabilizer of a point $x_0$. Each element of $H$ commutes with $S_{x_0}$, so
$h\coloneqq S_{x_0}$ lies in $Z(H)$. Throughout, we assume that the conjugates of $h$ generate
$G$, as they do in this presentation.

The group pair $(G,H)$ supplies quandle cohomology classes on $X$ via the Inoue--Kabaya (IK)
map. Let $A_X$ be the realization of the simplicial set with $n$-simplices $X^{n+1}/G$, where
$G$ acts diagonally and the faces and degeneracies delete and repeat coordinates, respectively.
For $n>0$, the group $H^n(A_X;B)$ is the Adamson relative cohomology of $(G,H)$ with constant
coefficients in $B$. The IK chain map~\cite[\S3.2]{IK} induces
\[
 \Phi_B\colon H^n(A_X;B)\longrightarrow H_Q^n(X;B).
\]
If $H$ is finite and multiplication by $|H|$ is an automorphism of $B$, the domain is naturally
isomorphic to $H^n(BG,BH;B)$ (Lemma~\ref{good:marked-pinch}).

We study products of classes in the image of $\Phi_B$ using the rack space $BX$. Its cubes are
labeled by elements of $X$, one for each coordinate direction, and its face identifications are
given by the quandle operation~\cite{FRS}. Quandle cohomology embeds naturally in the cohomology
of $BX$. The cup product of two positive-degree rack classes splits into two half-products,
according to which tensor factor of the cubical diagonal contains the edge in the first
coordinate~\cite{Clauwens,CFLM}.

The loop space of the double mapping cylinder $Y$ in \cite[Theorem 2.1]{Model} is homotopy
equivalent to the covering rack monoid of $X$. The two end inclusions, $BG\to Y$ and $BH\to Y$, have homotopy fibers modeled by $BX$ and $B(X;X)$, respectively
(Theorem~\ref{thm:model}). Each cube of the shadow rack space $B(X;X)$ carries an additional
label in $X$ at the initial vertex. The orbit map $BG\to A_X$ extends to $Y$ by a constant map
on the cylinder and on $BH$, and the extension induces a map $B(X;X)\to\Omega A_X$. On chains,
its adjoint $\Sigma B(X;X)\to A_X$ represents the reduced shadow IK map up to chain homotopy
over $\Z$ (Theorem~\ref{thm:comparison}).

The map $S^1\times B(X;X)\to BX$ of Fenn, Rourke, and Sanderson yields a continuous map
$\Sigma BX\to\Sigma A_X$ whose pullback corresponds, under the suspension isomorphisms, to the
composite of $\Phi_B$ with the inclusion into rack cohomology. Thus $\Phi_B$ commutes in
positive degrees with stable cohomology operations with constant coefficients
(Proposition~\ref{ik:stable}). If $H$ is finite and $|H|$ is invertible in the coefficient
ring, the image of the ordinary IK map on homology in degree $n\ge2$ is the image of the
homology suspension for $A_X$ (Proposition~\ref{good:image}).

Now suppose that $X$ and $G$ are finite, let $p$ be a prime with $p\nmid |H|$, and put
$R\coloneqq\Z_{(p)}$. We compare $BX$ with a James reduced product, a space of words in which
the basepoint is identified with the unit~\cite[\S1]{James}. The comparison is a zigzag of
continuous maps and induces a ring isomorphism on cohomology with coefficients in any
commutative $R$-algebra (Proposition~\ref{james:space}). Over a field of characteristic $p$,
the ordinary rack cohomology ring is noncanonically isomorphic to a shuffle algebra
\cite[Proposition~5.9]{Model}. We show that the comparison above carries the canonical
half-product on positive-degree rack cohomology to the half-shuffle on words
(Theorem~\ref{half:main}). This operation sums, with Koszul signs, the
interleavings that preserve the order in each word and begin with the first letter of the left
word. It follows that positive-degree rack cohomology is a free graded Zinbiel algebra. If a
mod-$p$ quandle class of positive even degree lifts modulo $p^2$, then the image of its
Pontryagin--Thomas value in rack cohomology is given by a word formula involving a divided
power (Corollary~\ref{prop:pontryagin-halfpower}).

With coefficients in a commutative ring, $\Phi_B$ annihilates cup products of positive-degree
classes (Corollary~\ref{lem:cupzero}). For the eight-element quandle
$X=\F_3^2\setminus\{0\}$ with $v*w=v+\det(v,w)w$ and coefficients $B=\Z_{(2)}$, the map
$\Phi_B\colon H^4(A_X;B)\to H_Q^4(X;B)$ is an isomorphism between cyclic groups of order $8$
(Theorem~\ref{binary:main}). If $a$ generates the domain, then $\Phi_B(a^2)=0$ and
$\Phi_B(a)^2$ has additive order $4$. Taking the Pontryagin square of $a\bmod 2$ and then
applying the IK map gives zero. Applying the IK map before taking the Pontryagin square gives
a nonzero class.

More generally, under a vanishing hypothesis on the low-degree $R$-homology of $A_X$, we
compute the additive orders of the cup powers of the IK image of a generator of a cyclic direct
summand in even-degree cohomology. We also compute the orders of the Pontryagin--Thomas values
of its reductions modulo powers of $p$ (Theorem~\ref{james:all-powers} and
Corollary~\ref{james:entire-tower}).

Section~\ref{sec:binary} applies these results to symplectic quandles over finite fields,
including the eight-element example. For each prime $p$, the nilpotence indices of degree-four
IK images over $\Z_{(p)}$ are unbounded in a family of fixed type, where the type is the common
order of the right translations (Corollary~\ref{symplectic:unbounded}). A family of cyclic
Alexander quandles of type $3$ exhibits the same phenomenon for degree-six IK images over
$\Z_{(7)}$ (Example~\ref{cyclic:sixth-degree}).

\section{Chain complexes and the marked model}\label{sec:ik-setup}
In this section, we review that the ordinary and shadow rack spaces model the homotopy fibers of the two end inclusions $BG\to Y$
and $BH\to Y$, respectively. Extending the orbit map by the constant map on the cylinder and $BH$ gives
$F:Y\to A_X$, which induces $B(X;X)\to\Omega A_X$ on the second fiber.

\subsection{Augmentations and chain complexes}
We work in the category of compactly generated weak Hausdorff spaces.

A quandle $X$ is a set with an idempotent, right self-distributive operation $x*y$ for which
every right translation $S_y:x\mapsto x*y$ is a bijection. Let $\Inn(X)$ denote the group
generated by the right translations. The quandle $X$ is connected if this group acts
transitively. The type of $X$, when it exists, is the least positive integer $t$ such that
$S_y^t=\id_X$ for all $y$. For permutations acting on the right, the product $gh$ applies $g$
first and then $h$, so that $x^{gh}=(x^g)^h$. An augmentation consists of a right action of a
group $G$ on $X$ and a map $\eta:X\to G$ satisfying
\begin{equation}\label{eq:augmentation}
 \eta(x^g)=g^{-1}\eta(x)g,\qquad
 x*y=x^{\eta(y)},\qquad x^{\eta(x)}=x\,.
\end{equation}
Throughout, we assume that the action is transitive and that $\eta(X)$ generates $G$. The
canonical augmentation is given by $G=\Inn(X)$ and $\eta(x)=S_x$. Fix $x_0\in X$ and write
\[
 H\coloneqq \Stab_G(x_0),\qquad h\coloneqq \eta(x_0),\qquad i:H\hookrightarrow G\,.
\]
We have $h\in Z(H)$, and the map $Hg\mapsto x_0^g$ identifies $X$ with $H\backslash G$. In these
coordinates the quandle operation is
\[
 Hg*Hk=Hgk^{-1}hk\,.
\]
For a right $G$-set $U$, put $u*x\coloneqq u^{\eta(x)}$. Write
$\boldsymbol x\coloneqq (x_1,\ldots,x_n)$ and define the cubical face maps $d_i^0,d_i^1$
($1\le i\le n$) and the chain boundary by
\begin{equation}\label{eq:rack-differential}
 \begin{aligned}
 d_i^0(u;\boldsymbol x)&\coloneqq (u;x_1,\ldots,\widehat{x_i},\ldots,x_n),\\
 d_i^1(u;\boldsymbol x)&\coloneqq (u*x_i;x_1*x_i,\ldots,x_{i-1}*x_i,x_{i+1},\ldots,x_n),\\
 \partial&\coloneqq \sum_{i=1}^n(-1)^i(d_i^0-d_i^1)\,.
 \end{aligned}
\end{equation}
Let $B(U;X)$ denote the space obtained by assigning an $n$-cube to each
$(u;\boldsymbol x)\in U\times X^n$ and gluing along the face maps in
\eqref{eq:rack-differential}. For the trivial one-point $G$-set $U=\{*\}$, its cellular chain
complex is the integral rack complex $C_n^R(X)\coloneqq \Z[X^n]$, with differential
$\partial_R$. Let $C_*^D(X)$ be the degenerate subcomplex generated by tuples with two equal
adjacent entries, and put $C_*^Q(X)\coloneqq C_*^R(X)/C_*^D(X)$. By the Litherland--Nelson
splitting \cite[Theorem 2.1]{LN}, this quotient map admits a chain section.
For a constant coefficient module $B$, write $H_*^\bullet(X;B)$ and $H_\bullet^*(X;B)$ for the
homology of $C_*^\bullet(X)\otimes B$ and the cohomology of $\Hom_{\Z}(C_*^\bullet(X),B)$,
respectively, where $\bullet=R,D,Q$. We use the cohomological differential
$\delta f\coloneqq f\partial$.

For $U=X$ the resulting complex is the shadow rack complex $S_n\coloneqq \Z[X\times X^n]$, with
differential $\partial_S$. In a generator $(r;x_1,\ldots,x_n)$, we call $r$ the region label and
each $x_i$ an arc label. Reading the first coordinate as the region label yields the degree
shift $\rho(x_1,\ldots,x_n)\coloneqq (x_1;x_2,\ldots,x_n)$ of Fenn--Rourke--Sanderson
\cite[Theorem 5.14]{FRS}, which gives $C_n^R(X)\cong S_{n-1}$ for $n\ge1$. Setting
$\rho|_{C_0^R(X)}\coloneqq 0$ gives $\partial_S\rho=-\rho\partial_R$.
In the shadow quandle normalization, tuples with two equal adjacent arc labels are degenerate.
Under $\rho$, the ordinary degenerate subcomplex also includes tuples whose region label equals
the first arc label.

Let $E(X)_\bullet$ be the simplicial set with $E(X)_n\coloneqq X^{n+1}$. For $0\le i\le n$,
define the face maps $d_i$ ($n\ge1$) and degeneracy maps $s_i$ ($n\ge0$) by
\[
 \begin{aligned}
 d_i(z_0,\ldots,z_n)&\coloneqq(z_0,\ldots,\widehat{z_i},\ldots,z_n),\\
 s_i(z_0,\ldots,z_n)&\coloneqq(z_0,\ldots,z_i,z_i,\ldots,z_n)\,.
 \end{aligned}
\]
This is the nerve of the pair groupoid on $X$, so its realization $E\coloneqq |E(X)_\bullet|$ is
contractible. Define the orbit space for the diagonal $G$-action and its normalized chain
complex by
\[
 A_X\coloneqq |E(X)_\bullet/G|,\qquad
 D_*(X)\coloneqq N_*(E(X)_\bullet/G;\Z)\,.
\]
Since $\eta(X)$ generates $G$ and $\eta(y)$ acts as $S_y$, the image of $G$ in the permutation
group of $X$ is $\Inn(X)$. Thus the $G$-orbits and the $\Inn(X)$-orbits agree in each simplicial
degree, and $A_X$ depends only on $X$. We write $q$ for both the quotient map $E\to A_X$ and its
induced chain map. Let $\varepsilon$ denote the augmentation sending every vertex in $S_0$ and
$D_0(X)$ to $1$, and let $\widetilde S$ and $\widetilde D$ be the corresponding reduced
complexes. We call the homology of $D_*(X)$ Adamson homology, as distinct from the Takasu
relative group homology $H_*(BG,BH)$ \cite{AdamsonComparison}.

Transitivity gives $A_X$ one vertex. Its edge loops $e_g=[x_0,x_0^{g^{-1}}]$ generate
$\pi_1(A_X)$, and the triangle $(x_0,x_0^{g^{-1}},x_0^{(gk)^{-1}})$ gives $e_ge_k=e_{gk}$.
Thus $g\mapsto e_g$ is a surjection. It kills $H$ since $H$ fixes $x_0$. As $h\in H$ and its conjugates
constitute the generating set $\eta(X)$, we obtain $\pi_1(A_X)=1$.

\subsection{A delooping under \texorpdfstring{$BG$}{BG}}

Write $BX\coloneqq B(\{*\};X)$ for the ordinary rack space. The Ishikawa--Tanaka quandle space
$B^QX$ is obtained from $BX$ by the following identifications in each cube with $x_i=x_{i+1}$.
Two points are identified if their coordinates outside positions $i,i+1$ agree and their sums
$t_i+t_{i+1}$ agree \cite[\S6.1]{IT}. The quotient map $\pi_X:BX\to B^QX$ maps nondegenerate
cells to the corresponding cells and sends degenerate cells to lower skeleta, so it induces
normalization on cellular chains. Thus $H^*(B^QX;B)=H_Q^*(X;B)$. The integral chain splitting
above implies that $\pi_X^*$ is injective for every constant coefficient module. For
coefficients in a commutative ring it preserves the ordinary cup product.

We write the covering rack space and the shadow space as
\[
 M\coloneqq B(G;X),\qquad \mathcal N\coloneqq B(X;X)=H\backslash M\,.
\]
The quotient $v:M\to\mathcal N$ is the covering map given by the left $H$-action. The space $M$
is a connected topological monoid containing $G$ as its group of vertices
\cite{Clauwens,Mostovoy,Model}. Define two group homomorphisms by
\[
 a:H\times\Z\to G,\quad a(u,m)\coloneqq i(u),\qquad
 b_h:H\times\Z\to H,\quad b_h(u,m)\coloneqq uh^m\,.
\]
They give the marked double mapping cylinder
\begin{equation}\label{eq:Y}
 Y\coloneqq \hocolim\bigl[BG\xleftarrow{Ba}B(H\times\Z)\xrightarrow{Bb_h}BH\bigr]\,.
\end{equation}
The marking is the inclusion $j:BG\to Y$ of the left end.

\begin{theorem}[Marked integral model]\label{thm:model}
The space $Y$ is simply connected, and there is a homotopy equivalence $BM\simeq Y$ under
$BG$. Moreover,
\[
 BX\simeq\hofib(BG\to Y),\qquad M\simeq\Omega Y,\qquad
 \mathcal N\simeq\hofib(BH\to Y),
\]
where $BG\to Y$ and $BH\to Y$ are the inclusions of the two ends.
\end{theorem}
The assertions concerning $Y$, $BM$, $BX$, and $M$ are proved in \cite[Theorem 2.1]{Model}. For
$\mathcal N$, the two-sided bar construction gives
$H\backslash M\simeq B(*,H,M)\simeq\hofib(BH\to BM)$. Under the comparison $BM\simeq Y$, the
cylinder homotopy along $u\mapsto(u,0)$ identifies the map from $BH$ in this construction
with the inclusion of the right end. We use the Mostovoy comparison and the equivalence
specified in that reference by the edge $w(t)\coloneqq[1;t,x_0]$, with the comparison at the
$BH$ end and along the cylinder fixed as in Appendix~\ref{appcomp:section}.

In the homogeneous bar model of $BG$, define the orbit map by
\[
 q_G:BG\longrightarrow A_X,\qquad
 [g_0,\ldots,g_n]\longmapsto[x_0^{g_0},\ldots,x_0^{g_n}]\,.
\]
Every vertex of $EH$ maps to $x_0$, so $q_G|_{BH}$ and $q_G\circ Ba$ are constant. Extending
$q_G$ by the constant map on the cylinder and the right end $BH$ gives a continuous map
\begin{equation*}
 F:Y\longrightarrow A_X.
\end{equation*}
Take the vertex of $BG$ as the basepoint of $Y$, and let $j_H:BH\hookrightarrow Y$ denote
the inclusion of the right end. Define $\kappa$ by the composite
\[
 \mathcal N\xrightarrow{\ \simeq\ }\hofib(j_H)
 \xrightarrow{\ (b,\gamma)\,\mapsto\, F\circ\gamma\ }\Omega A_X,
\]
using the marked equivalence of Theorem~\ref{thm:model}. Here $b\in BH$ and $\gamma$ is a path
from the basepoint of $Y$ to $j_H(b)$. Both endpoints map to the vertex of $A_X$, and the
cylinder path representing the basepoint of the homotopy fiber maps to the constant loop.
Write $\widehat\kappa:\Sigma\mathcal N\to A_X$ for the adjoint.

\section{Comparison with the original IK map in all degrees}\label{sec:ik-comparison}

We realize the reduced shadow IK map up to integral chain homotopy by
$\widehat\kappa:\Sigma\mathcal N\to A_X$. Products of two positive-degree classes on $\Sigma\mathcal N$
vanish, so $\Phi_B(a\smile b)=0$ for $a,b\in H^{>0}(A_X;B)$ over any commutative ring $B$.

For $1\le i\le n$, define $R_i:X^{n+1}\to X^{n+1}$ by
\[
 (r,x_1,\ldots,x_n)\longmapsto
 (r*x_i,x_1*x_i,\ldots,x_{i-1}*x_i,x_i,\ldots,x_n)\,.
\]
The self-distributive law gives $R_iR_j=R_jR_i$. Define the shadow finite-difference map
$\Psi_n:S_n\to C_n^{\mathrm{cell}}(E;\Z)$ by the following formula, with the empty product in
degree zero giving the identity on vertices.
\begin{equation}\label{eq:Psi}
 \Psi_n(r;\boldsymbol x)\coloneqq \prod_{i=1}^n(1-R_i)(r,x_1,\ldots,x_n)\,.
\end{equation}
Let $h_p$ be the cone operator that prepends a vertex $p\in X$. Define the shadow IK map by
\begin{equation}\label{eq:Theta}
 \Theta_n\coloneqq q h_p\Psi_n:S_n\longrightarrow D_{n+1}(X)\,.
\end{equation}
For the ordinary version, read $x_1$ as the region label and set, for $n\ge1$,
\begin{equation}\label{eq:ordinaryIK}
 \varphi_n(x_1,\ldots,x_n)
 \coloneqq \Theta_{n-1}(x_1;x_2,\ldots,x_n)
 =qh_p\Psi_{n-1}(x_1;x_2,\ldots,x_n)\,.
\end{equation}
Set $\varphi_0\coloneqq 0$. We refer to this as the original Inoue--Kabaya map \cite[\S3.2]{IK}.

\begin{remark}[Comparison with Inoue's normalized map]\label{rem:inoue-specialization}
Inoue \cite[Theorem 4.1]{Inoue2026} constructs a chain map $\psi$ for the $\Z$-shadow, on which
each arc label acts by $u\mapsto u+1$, and shows that it descends to the quandle normalization.
Put $\phi(g)\coloneqq h^{-1}gh$. Then $\phi|_H=\id_H$, and the operation $H\phi(gk^{-1})k$ in
that reference agrees with $Hgk^{-1}hk$. The map $Hg\mapsto g^{-1}H$ reconciles the left and
right conventions for the Adamson complex. For $x_j=Hg_j$ and $m\in\Z$, we have
$(\phi^{-m}(g_j))^{-1}H=h^m g_j^{-1}H$, so the simultaneous left action of $h^m$ disappears in
the orbit quotient. With the same auxiliary vertex in both constructions, the image of
$\psi_n(u;x_1,\ldots,x_n)$ in the normalized simplicial complex is $-\varphi_n(x_1,\ldots,x_n)$
for $n\ge1$. In degree zero the $\Z$-shadow map satisfies $\psi_0(u)=-(H)$. The ordinary complex
in that reference is zero in degree zero. Here we take $C_0^R(X)=\Z$ and $\varphi_0=0$.

In the proof of Theorem~4.2 of that reference, on independence of the auxiliary vertex, the
prism identity has an additional edge term in degree one.
Lemma~\ref{lem:finite-difference} shows how to account for this term in the ordinary IK map.
\end{remark}

\begin{lemma}\label{lem:finite-difference}
For every $n\ge1$ the following identities hold.
\begin{gather*}
 \partial\Psi_n=\Psi_{n-1}\partial_S,\qquad
 q\Psi_n=0,\\
 \partial\Theta_n=-\Theta_{n-1}\partial_S,
 \qquad \partial\varphi_n=\varphi_{n-1}\partial_R\,.
\end{gather*}
The map $\varphi$ vanishes on $C_*^D(X)$ and therefore descends to $C_*^Q(X)$. For both the
reduced shadow map and the ordinary IK map, different choices of $p$ give chain-homotopic
maps.
\end{lemma}
\begin{proof}
Equation~\eqref{appcomp:Psi} and Lemma~\ref{appcomp:orbit-cancellation} give the first line. The
identities $\partial h_p+h_p\partial=\id-[p]\eps$ and $q\Psi_{>0}=0$ give
$\partial\Theta_n=-\Theta_{n-1}\partial_S$ for $n\ge1$. In degree zero, $\partial\Theta_0=0$
because the orbit space has a single vertex. The identity $\partial_S\rho=-\rho\partial_R$
cancels the minus sign for the ordinary map.

Suppose that two adjacent entries of $(r,x_1,\ldots,x_n)$ are equal. Each $R_j$ translates
both entries or neither unless the pair straddles the cut immediately before $x_j$. In that
case equality follows from $x_j*x_j=x_j$. Thus the two entries agree in every summand of
\eqref{eq:Psi} and remain equal after $h_p$ prepends $p$. This includes $r=x_1$, so \eqref{eq:ordinaryIK}
kills all ordinary degeneracies.

Put $K\coloneqq h_ph_{p'}$. Then $\partial K-K\partial=h_{p'}-h_p+h_p[p']\eps$. For the reduced
shadow map, $qK\Psi$ gives the homotopy. For the ordinary map, writing the auxiliary vertex in
parentheses, we have
\[
\partial(qK\Psi\rho)+(qK\Psi\rho)\partial_R
=\varphi(p')-\varphi(p)+q(p,p')\eps\rho.
\]
The last term occurs only in degree one. Since $H_1(A_X;\Z)=0$, choose $b\in D_2(X)$ with
$\partial b=q(p,p')$ and subtract $b\eps\rho$ in degree one. The equality $\eps\rho\partial_R=0$
ensures the degree-two identity as well. Prepending vertices preserves equality of adjacent
entries as above, so the homotopy also descends to the normalization.
\end{proof}

We henceforth write $\varphi^R$ for the ordinary rack chain map $\varphi$ and $\varphi^Q$ for
its descent to the quandle normalization. The induced cohomology homomorphism with constant
coefficients $B$ is $\Phi_B$. We use the convention $d(\susp c)\coloneqq -\susp(dc)$ for
suspension of complexes and regard $\Theta$ as a chain map $\susp\widetilde S\to\widetilde D$.

\begin{theorem}[Comparison in all degrees]\label{thm:comparison}
For every connected discrete quandle with a transitive augmentation whose image generates $G$,
the homomorphism induced by \eqref{eq:Theta}
\[
 \Theta_{n*}:\widetilde H_n(\mathcal N;\Z)
 \longrightarrow\widetilde H_{n+1}(A_X;\Z)
\]
is the composite of the suspension isomorphism with $\widehat\kappa_*$. The comparison is
realized by a chain homotopy over $\Z$ on free chain complexes and is natural with respect to
homomorphisms of constant coefficient modules.

In degrees $n\ge2$, the ordinary version is obtained by the degree shift
$H_n^R(X)\cong H_{n-1}(\mathcal N)$ and descends to the quandle normalization. Its induced map
in degree one is zero because $H_1(A_X)=0$. Proposition~\ref{ik:stable} gives a realization of
this degree shift after suspension.
\end{theorem}
\begin{proof}
In the model $Q$ of Appendix~\ref{appcomp:section}, the cylinder path between the basepoints
of $BG$ and $BH$ becomes constant. Write $j_H:BH\to Q$ for the map from the $BH$ end in this
model. Then $\widehat\kappa$ is represented by
$\Sigma\hofib(j_H)\to\Cofib(j_H)\to A_X$.
Here $\Cofib(j_H)$ is the reduced mapping cone, with the cone on the basepoint of $BH$
collapsed. The adjoint of the first map takes a fiber path over $b\in BH$ and follows it
by the path from $j_H(b)$ to the cone point, forming a based loop in the reduced cone.
On the model $B(*,H,M)\simeq\hofib(j_H)$ of Appendix~\ref{appcomp:section}, the first map
is represented on suspended chains by $I$ in \eqref{appcomp:I}.
The second map extends $F$ in the model $Q$ by the constant map on the cone.
Let $K$ be the chain homotopy of \eqref{appcomp:beta-J} transported to $Q$ under $BH$.
The extension $\overline K(t,\susp u)=K(t)$ satisfies the cone homotopy identity
because $Kj_Hu=0$ on reduced chains by relativity to $BH$.
Likewise $\beta j_H=0$ on reduced chains, so $\overline\beta(t,\susp u)=\beta(t)$
represents the second map.

Take $p=x_0$. By $\beta i=q$ in \eqref{appcomp:beta} and $\beta\mathsf P=0$ in
\eqref{appcomp:prism-homotopy}, for $w\in\widetilde W$ we have
\[
 \overline\beta I(\susp w)=\beta(C-\mathsf P)w=q h_{x_0}\Psi\rho_H(w)=\Theta(\susp\rho_H(w))\,.
\]
The comparison map $\rho_H$ for the covering quotient preserves the augmentations and remains a
quasi-isomorphism of bounded-below free integral chain complexes after reduction. By
Lemma~\ref{good:coefficient-change} with $R=\Z$, it is an integral chain homotopy equivalence.
The comparison is therefore natural in the constant coefficients. The ordinary version follows
from \eqref{eq:ordinaryIK}, and normalization and change of auxiliary vertex follow from
Lemma~\ref{lem:finite-difference}.
\end{proof}

\begin{corollary}\label{lem:cupzero}
For any commutative ring $B$ of constant coefficients and any positive-degree classes
$a,b\in H^{>0}(A_X;B)$, we have $\Phi_B(a\smile b)=0$.
\end{corollary}
\begin{proof}
By Theorem~\ref{thm:comparison} and \eqref{eq:ordinaryIK}, the class $\pi_X^*\Phi_B(a\smile b)$
is the degree shift of $\widehat\kappa^*(a\smile b)$. The latter is the product
$\widehat\kappa^*a\smile\widehat\kappa^*b$ of positive-degree classes on a suspension and hence
vanishes. The claim follows from the injectivity of $\pi_X^*$.
\end{proof}
Products of classes in the image are treated in Section~\ref{james:section}.

\section{The IK map after inverting the stabilizer order}\label{sec:ik-good}
In this section, we show that, over a localization $R$ of $\Z$ inverting the finite stabilizer order $|H|$, the IK image in
degrees $n\ge2$ equals the image of homology suspension for $A_X$. If $H_j(A_X;R)=0$ for
$2\le j<n$, the quandle IK map in degree $n$ is an isomorphism for every $R$-module of coefficients.

All spaces are connected and have CW type. A subscript $R$ denotes the
$R$-localization of a simply connected space. For the orbit space we write
$A_{X,R}\coloneqq(A_X)_R$.

\begin{lemma}[Change of coefficients]\label{good:coefficient-change}
Let $R$ be any localization of $\Z$. A quasi-isomorphism between bounded-below complexes of free
$R$-modules is a chain homotopy equivalence. In particular, an $R$-homology equivalence between
spaces of CW type induces isomorphisms in homology and cohomology with coefficients in any
constant $R$-module $B$.
\end{lemma}
\begin{proof}
The mapping cone is a bounded-below acyclic complex of free $R$-modules. Since $R$ is a PID, its
cycle modules are free, and degreewise splittings give a chain contraction. Tensor products and
$\Hom_R(-,B)$ preserve chain homotopy equivalences.
\end{proof}

\begin{lemma}[The localized splitting and the pinch map]\label{good:marked-pinch}
Write $L\coloneqq \Cofib(BH\to BG)$ for the homotopy cofiber of the inclusion $BH\to BG$. Both
the orbit map $q_G:BG\to A_X$ and the map $\lambda:L\to A_X$ obtained from its specified
constant restriction on $BH$ induce isomorphisms in $R$-homology. The comparison respecting the specified maps on $BH$ and the cylinder gives
\[
 L_R\simeq A_{X,R},\qquad
 Y_R\simeq(S^2\vee A_X)_R,\qquad
 F_R\simeq(\operatorname{pr}_{A_X})_R\,.
\]
The last map is the localization of the pinch map that collapses the sphere summand to a point.
\end{lemma}
\begin{proof}
Filter the Borel model $B(E,G,*)$ by simplicial dimension in $E$. Since each simplex stabilizer
fixes the simplex pointwise, the corresponding $E^1_{s,t}$ is the direct sum, over orbits of
$s$-simplices, of the degree-$t$ homology of their stabilizers \cite[Chapter VII, \S7,
(7.7)]{Brown}. Each stabilizer is a subgroup of a conjugate of $H$, so its order is invertible
in $R$. Therefore the rows with $t>0$ vanish and the row with $t=0$ is the chain complex of
$A_X$. Hence $B(E,G,*)\to A_X$ induces an $R$-homology isomorphism, and the Borel model
comparison in Appendix~\ref{appcomp:section} gives the same conclusion for $q_G$. By \cite[Lemma
3.1]{Model}, the quotient $\ell:BG\to L=BG/BH$ is an $R$-homology equivalence and $L$ is simply
connected. Since $\lambda\ell=q_G$, the map $\lambda$ is also an $R$-homology equivalence. As
$A_X$ is simply connected, $\lambda_R$ is a homotopy equivalence.

The map of spans in the proof of \cite[Theorem 3.2]{Model}, from
$[BG\leftarrow B(H\times\Z)\to BH]$ to $[L\leftarrow B\Z\to *]$, has components $\ell$,
$B\operatorname{pr}_2$, and the collapse of $BH$. Both squares commute strictly because
$L=BG/BH$, and each component is an $R$-homology equivalence since
$\widetilde H_*(BH;R)=0$.
Collapse the cylinder over the basepoint and identify $\Sigma B\Z\simeq S^2$, using $1\in\Z$ to
specify the positive orientation. Composing the induced map with $\lambda$ on the cofiber
summand and interchanging the wedge factors gives an $R$-homology equivalence
$g:Y\to S^2\vee A_X$. On the $BG$ end we have $gj=i_{A_X}q_G$, where $i_{A_X}$ includes the
$A_X$ summand. The $BH$
end maps to the basepoint and the cylinder maps to the sphere. Thus
$\operatorname{pr}_{A_X} g=F$ with the specified constant restrictions. Localization gives the
result.
\end{proof}

The covering quotient $v$ induces the following isomorphism for every $R$-module $\Lambda$.
\begin{equation}\label{good:cover-transfer}
 v_*:H_*(M;\Lambda)\xrightarrow{\cong}H_*(\mathcal N;\Lambda)\,.
\end{equation}
Indeed, the deck transformation $m\mapsto um$ is homotopic to the identity by left
multiplication along a path in $M$ from $1$ to $u$. Consequently both composites of $v_*$ with
the transfer are multiplication by $|H|$, which is invertible on $\Lambda$ \cite[Propositions
22--23]{Clauwens}~\cite[Proposition 5.1]{Model}.

\begin{proposition}[Image in all degrees]\label{good:image}
Let $\Lambda$ be an $R$-module and let $n\ge2$. Let $c_{M,R}:M\to\Omega Y_R$ be the composite of
the marked bar--loop comparison $c_M:M\to\Omega Y$ of Appendix~\ref{appcomp:section} with
$\Omega(Y\to Y_R)$. The rack version of the IK map satisfies
\begin{equation}\label{good:signed-ik}
 (A_X\to A_{X,R})_*\circ\varphi^R_{n*}
 =-\sigma_{A_{X,R}}\circ(\Omega F_R)_*\circ(c_{M,R})_*\circ v_*^{-1}\circ\rho_*\,.
\end{equation}
Here $\sigma$ denotes homology suspension. The rack and quandle IK maps both have image equal to
that of suspension.
\[
 \operatorname{im}\varphi^R_{n*}
 =\operatorname{im}\varphi^Q_{n*}
 =\operatorname{im}\bigl(
 \sigma_{A_X}:H_{n-1}(\Omega A_X;\Lambda)\longrightarrow H_n(A_X;\Lambda)
 \bigr)\,.
\]
\end{proposition}
\begin{proof}
By Theorem~\ref{thm:comparison}, $\varphi^R_{n*}=\sigma_{A_X}\kappa_*\rho_*$. Substituting the formula for
$\kappa v$ at the end of Appendix~\ref{appcomp:section} and localizing gives \eqref{good:signed-ik}.
The minus sign comes from reversing the suspension coordinate.
In \eqref{good:signed-ik}, $\rho_*$ and $v_*$ are isomorphisms by the degree shift in
Section~\ref{sec:ik-setup} and \eqref{good:cover-transfer}. Since $c_M$ is an equivalence and
localization commutes with looping simply connected spaces \cite[proof of Theorem 3.2]{Model},
$(c_{M,R})_*$ is also an isomorphism, and localization identifies the two suspension images
with $\Lambda$ coefficients by Lemma~\ref{good:coefficient-change}. By Lemma~\ref{good:marked-pinch}, $F_R$ has
a homotopy section, so $(\Omega F_R)_*$ is surjective. The rack assertion follows, and the
quandle assertion follows by Lemma~\ref{lem:finite-difference} and the split projection
$C_*^R\to C_*^Q$.
\end{proof}

\begin{theorem}[Detection in the first possible nonzero degree]\label{good:degree3}
Under the above assumptions, let $n\ge2$ and suppose that $H_j(A_X;R)=0$ for $2\le j<n$. Then
for every $R$-module $\Lambda$,
\[
 H_j^Q(X;\Lambda)=0\quad(2\le j<n),\qquad
 \varphi_{n*}:H_n^Q(X;\Lambda)\xrightarrow{\ \cong\ }H_n(A_X;\Lambda)\,.
\]
The isomorphism is induced by the IK map fixed in Theorem~\ref{thm:comparison}.
\end{theorem}
\begin{proof}
First take coefficients in $R$. The space $A_{X,R}$ is $(n-1)$-connected. Choose a CW model for
$A_{X,R}$ whose positive-dimensional cells have dimension at least $n$, and put
$W\coloneqq S^2\vee A_{X,R}$. The map $S^2\vee A_X\to W$ is an $R$-homology equivalence. By
Lemma~\ref{good:marked-pinch}, under $Y_R\simeq W_R$ the map $F_R$ corresponds to the
localization of the projection onto the second summand. The relative cells of the inclusion $W\hookrightarrow S^2\times A_{X,R}$ have dimension at least $n+2$, so
$\Omega W\to\Omega S^2\times\Omega A_{X,R}$ is $n$-connected. The one-element rack space
$B\{*\}\simeq\Omega S^2$ has one cell in each dimension and zero cellular differential
\cite[Theorem 5.12]{FRS}~\cite{James}. Using this local splitting, \eqref{good:cover-transfer},
the rack--shadow degree shift, and the K\"unneth formula in low degrees, we obtain
\[
 H_j^R(X;R)\cong R\quad(1\le j<n),\qquad
 H_n^R(X;R)\cong R\oplus H_{n-1}(\Omega A_{X,R};R)\,.
\]
By \eqref{good:signed-ik}, the rack IK map in degree $n$ is projection onto the second summand
followed by $-\sigma_{A_{X,R}}$. The Hurewicz isomorphisms in the first possible nonzero degrees
identify this suspension with an isomorphism $H_{n-1}(\Omega A_{X,R};R)\cong H_n(A_{X,R};R)$.
Identifying the target with $H_n(A_X;R)$ via localization, we see that $\varphi_{n*}^R$ is
surjective with kernel isomorphic to $R$.

The map $BX\to B\{*\}$ that forgets the labels sends the degenerate cycle
$d_j\coloneqq(x_0,\ldots,x_0)$, for $j\ge2$, to the generator in degree $j$. Hence $d_j$
generates the rank-one module $H_j^R(X;R)$ for $2\le j<n$. In degree $n$, the forgetful map
restricted to the rank-one kernel $\ker\varphi_{n*}^R$ sends $d_n$ to a generator, so $d_n$
generates the kernel. Thus the image of $H_j^D$ is all of $H_j^R$ for $2\le j<n$ and is the kernel of the rack IK
map for $j=n$. The integral chain splitting gives the result over $R$. Connectedness also gives $H_1^Q(X;R)=R$, and
$H_1(A_X;R)=0$. The universal coefficient $\Tor_R$ terms therefore vanish for $2\le j\le n$.
Naturality of the universal coefficient theorem gives the result for every $\Lambda$.
\end{proof}

\section{Comparison with the James reduced product}\label{sec:james-comparison}

For finite $X,G$ and $p\nmid |H|$, we will compare the rack space with a James reduced product
through a common homotopy fiber. The maps identify their cohomology rings over any commutative
$\Z_{(p)}$-algebra.

Put $R\coloneqq\Z_{(p)}$. Using
$A_{X,R}=(A_X)_{(p)}$, set
\[
 \mathcal{K}\coloneqq \Omega A_{X,R},\qquad W\coloneqq S^2\vee A_{X,R},
 \qquad i:A_{X,R}\hookrightarrow W,\qquad V\coloneqq \Sigma \mathcal{K}_+\,.
\]
Here $\mathcal{K}_+\coloneqq \mathcal{K}\sqcup\{*\}$ is obtained by adjoining a disjoint
basepoint to $\mathcal K$. We use CW models with the basepoint a $0$-cell. The James reduced
product $J(V)$ is the free topological monoid on $V$ with the basepoint of $V$ identified with
the unit \cite[Theorem (1.11)]{James}.

\begin{proposition}[Comparison through a common homotopy fiber]\label{james:space}
There is a zigzag of continuous maps
\begin{equation}\label{james:zigzag}
 BX\xrightarrow{\ \beta\ }F_i\xleftarrow{\ \alpha\ }J(V),
 \qquad F_i\coloneqq \hofib(i)
\end{equation}
in which $\beta$ is an $R$-homology equivalence and $\alpha$ is a weak equivalence. Both maps induce isomorphisms in homology and cohomology with coefficients in any constant $R$-module $\Lambda$. Write
\begin{equation}\label{half:coefficient-comparison}
 \mathcal T_\Lambda\coloneqq \alpha^*(\beta^*)^{-1}:
 H_R^*(X;\Lambda)\xrightarrow{\cong}H^*(J(V);\Lambda)
\end{equation}
for the resulting cohomology isomorphism. The comparison is natural in $\Lambda$ and intertwines
natural cohomology operations whose source and target coefficients are constant $R$-modules.
It preserves cup products for every commutative $R$-algebra $\Lambda$. Formula~\eqref{half:alpha} defines $\alpha$ using right conjugation
and loop reversal.
\end{proposition}
\begin{proof}
Let $\ell_p:A_X\to A_{X,R}$ be the localization map, and use the map $g$ in the proof of
Lemma~\ref{good:marked-pinch} to define $\chi\coloneqq (\id_{S^2}\vee\ell_p)g:Y\to W$. Put
$\bar q_G\coloneqq \ell_pq_G$, so that $\chi j=i\bar q_G$. For the equivalences
$BX\simeq B(*,G,M)\simeq\hofib(j)$ from the proof of
\cite[Theorem 2.1]{Model}, use the prism comparison of Appendix~\ref{appcomp:section} with
action space $*$ and subgroup $G$. On $M$, this comparison uses the reverse of the bar edge
$\ell_a$ as the path from the basepoint to the point on the $BG$ end. Composing the
equivalences with the map of homotopy fibers
$(b,\gamma)\mapsto(\bar q_G(b),\chi\circ\gamma)$ induced by the square $\chi j=i\bar q_G$
gives $\beta$. The map $\bar q_G$ on the bases and $\Omega\chi$ on the fibers are $R$-homology
equivalences, by Lemma~\ref{good:marked-pinch} and the compatibility of localization with
looping for simply connected spaces. The two fibrations are pullbacks of the path fibrations of
the simply connected spaces $Y$ and $W$, so their coefficient systems are constant.
Lemma~\ref{good:coefficient-change} and comparison of the Serre spectral sequences give the
assertion about $\beta$.

Let $\mathcal{K}_\bullet$ be the Kan loop group of a reduced simplicial model of $A_{X,R}$, and
henceforth model $\mathcal{K}$ by its realization. We use the resulting based equivalence
$B\mathcal{K}\simeq A_{X,R}$ \cite[Ch.~V, Corollary~6.4]{GoerssJardine}. Let $S_\bullet$ be the
simplicial circle and let $\mathsf F$ denote Milnor's free simplicial group construction
\cite[\S1]{Milnor}. Put
\[
 P_\bullet\coloneqq \mathsf F(S_\bullet),\qquad
 \mathcal G_\bullet\coloneqq P_\bullet*\mathcal{K}_\bullet,\qquad
 \mathcal L_\bullet\coloneqq \ker(\mathcal G_\bullet\to \mathcal{K}_\bullet)
\]
and use the same letters for the realizations. The wedge description of classifying spaces of free products in each simplicial degree, together with $BP\simeq S^2$ \cite[Ch.~V, Corollary~6.17(2)]{GoerssJardine}, gives $B\mathcal G\simeq BP\vee B\mathcal{K}\simeq W$ compatibly with both inclusions. Under this equivalence $B\mathcal{K}\to B\mathcal G$ corresponds to $i:A_{X,R}\hookrightarrow W$. We use $\mathcal G$ as a model for $\Omega W$ via the standard bar--loop comparison \cite[Proposition~4.17(1), Theorem~5.6]{Vogt}. The normal form for a free product identifies the kernel by the following simplicial group isomorphism.
\begin{equation}\label{half:kernel}
 \mathsf F(S_\bullet\wedge \mathcal{K}_{\bullet,+})\xrightarrow{\cong}\mathcal L_\bullet,
 \qquad(t,k)\longmapsto k^{-1}t k\,.
\end{equation}
In each simplicial degree, send words in the James construction to positive words in Milnor's
free group, preserving the order of the letters \cite{James,Milnor}. Since
$V=|S_\bullet\wedge\mathcal K_{\bullet,+}|$ is connected, the realization of this inclusion,
$J(V)\to|\mathsf F(S_\bullet\wedge\mathcal K_{\bullet,+})|$, is a weak equivalence \cite[online
lecture notes, Ch.~3, Theorem~3.5]{WuSimplicial}. Composing with the realization of
\eqref{half:kernel} gives $\Gamma:J(V)\xrightarrow{\simeq}\mathcal L\hookrightarrow\mathcal G$.

Paths in $F_i$ run from the basepoint to $i(a)$. Let $\delta:\mathcal G\to F_i$ be the map that
uses a based loop as the path over the basepoint. The composite
$\mathcal{K}\to\mathcal G\to\mathcal{K}$ is the identity, and
$\mathcal G\cong\mathcal L\times \mathcal{K}$ as spaces. The homotopy exact sequence of
$\mathcal{K}\to\mathcal G\to F_i$ shows that $\delta|_{\mathcal L}$ induces isomorphisms on the
positive-degree homotopy groups. Since $\mathcal L$ and $F_i$ are connected, the map
\begin{equation}\label{half:alpha}
 \alpha\coloneqq \delta\circ\operatorname{inv}_{\mathcal G}\circ\Gamma
\end{equation}
is a weak equivalence. Here $\operatorname{inv}_{\mathcal G}$ is group inversion,
corresponding to reversal of the loop represented by the entire word. The isomorphisms with general coefficients follow from Lemma~\ref{good:coefficient-change}, and compatibility with operations and cup products follows from naturality.
\end{proof}

\section{The rack half-product and the James construction}\label{half:section}
In this section, we show that, for finite $X,G$ and $p\nmid |H|$, the comparison of Section~\ref{sec:james-comparison}
identifies the rack half-product in positive degrees with the James half-shuffle over any field
of characteristic $p$. On homology we recover the signed sum over pairs of nonempty subwords with the
first letter on the left.

\subsection{The cubical diagonal and the shadow half-diagonal}\label{sec:loop-structures}
Clauwens's monoid $M$ is a generalized James construction~\cite[\S2.6, Propositions
8--9]{Clauwens}. Let $E_x$ be the edge from the unit vertex to $\eta(x)$, let $T$ denote the
unital tensor algebra, and put $\mathcal A\coloneqq C_*^{\mathrm{cell}}(M;\Z)$. Orient each
cube by the order of its edge factors. As a graded module,
\begin{equation}\label{loop:algebra-normal-form}
 \mathcal A\cong\Z[G]\otimes_\Z T(\Z[X]),\qquad
 [g;x_1,\ldots,x_n]=gE_{x_1}\cdots E_{x_n}\,.
\end{equation}

\begin{proposition}[The chain bialgebra of the covering monoid]\label{loop:dg-bialgebra}
The differential graded algebra $\mathcal A$ is given by the group relations on vertices and the
formulas
\begin{equation}\label{loop:cross-relation}
 E_xg=gE_{x^g},\qquad
 d(g)=0,\qquad d(E_x)=\eta(x)-1\,.
\end{equation}
The cubical diagonal gives the coproduct
\begin{equation}\label{loop:dg-coproduct}
 \Delta(g)=g\otimes g,\qquad
 \Delta(E_x)=E_x\otimes\eta(x)+1\otimes E_x\,.
\end{equation}
With this coproduct and counit $\varepsilon(g)=1$, $\varepsilon(E_x)=0$, the complex
$\mathcal A$ is a differential graded bialgebra. To obtain $\mathcal A$ from the bialgebra of
Covez, Farinati, Lebed, and
Manchon~\cite[Theorem 1, Remark 6]{CFLM}, first invert its degree-zero generators. Then
impose the degree-zero relations induced by $\As(X)\to G$, $g_x\mapsto\eta(x)$. The associated
group $\As(X)$ has defining relations $g_{x*y}=g_y^{-1}g_xg_y$. The two conventions for edge
generators are related by $e_x\mapsto-E_x$.
\end{proposition}
\begin{proof}
Clauwens's product on a $1$-cell and a vertex gives $E_xg=gE_{x^g}$
\cite[Proposition 3, (22)]{Clauwens}. The boundary is zero on
vertices and sends $E_x$ to $\eta(x)-1$, completing \eqref{loop:cross-relation}.
The cubical Serre diagonal gives \eqref{loop:dg-coproduct}
\cite[Chap.~II, \S1, p.~441]{SerreCubical}. On products of cells,
the boundary satisfies the graded Leibniz rule, and the diagonal is multiplicative for the
graded tensor product. The identity $g\eta(x^g)=\eta(x)g$ gives compatibility with the cross
relation, and the bialgebra identities follow from the same formulas on generators.
\end{proof}

The antipode is given by $S(g)=g^{-1}$ and $S(E_x)=-E_x\eta(x)^{-1}$. The standard bar--loop
comparison of Theorem~\ref{thm:model} preserves multiplication up to homotopy \cite[Proposition
4.17(1), Theorem 5.6]{Vogt}. Thus on homology over a field the product corresponds to loop
concatenation and the antipode to loop reversal.

The ordinary rack complex is the complex of coinvariants for the left $G$-action.
\begin{equation*}
 C_*^R(X)\cong\Z\otimes_{\Z[G]}\mathcal A\,.
\end{equation*}
The diagonal descends to the quotient, and its dual gives the ordinary cup product \cite[Lemma
7, Remark 8]{CFLM}. The identity $[E_xg]=[E_{x^g}]$, together with $[g]=[1]$, shows that the
product of the covering monoid generally fails to descend.

Let $B$ be a commutative ring with unit. We use the unsigned tensor pairing
$(f\otimes g)(c\otimes d)=f(c)g(d)$ for homogeneous cochains $f,g$ and chains $c,d$ of matching
degrees $|f|=|c|$ and $|g|=|d|$. The graded tensor pairing used by Clauwens~\cite[(32),
(34)]{Clauwens} and in \cite[Corollary 17]{CFLM} includes a factor of $(-1)^{rs}$ when the
two degrees are $r$ and $s$. Clauwens uses the cochain differential $f\mapsto f\partial$
\cite[(79)]{Clauwens}. With the same labeling coordinates, the boundary in \cite[(2),
Proposition 10]{CFLM} is $\partial_{\mathrm{CFLM}}=-\partial$. Its graded dual sends a
degree-$q$ cochain $f$ to $(-1)^q f\circ\partial_{\mathrm{CFLM}}$ (see \cite[Example 18]{CFLM}
in degree one). Thus $\delta_{\mathrm{CFLM}}^q=(-1)^{q+1}\delta$ for our differential
$\delta f=f\partial$, and $\theta_q^+\coloneqq (-1)^{q(q+1)/2}\id$ satisfies
$\delta_{\mathrm{CFLM}}^q\theta_q^+=\theta_{q+1}^+\delta$. The sign in degree $r+s$ is
$(-1)^{rs}$ times the product of the signs in degrees $r,s$, so $\theta^+$ carries the product
and each half-product from the graded to the unsigned tensor pairing.

For positive-degree rack cochains $f,g$ with coefficients in $B$, define the canonical
half-product $f\prec g$ using the expansion of the cubical diagonal \eqref{loop:dg-coproduct}.
Retain the terms in which the edge in the first coordinate goes to the left factor and both
factors have positive degree. Keep their Koszul signs and use the unsigned pairing with $f\otimes g$ after passing to the
rack complex. On coinvariants this formula uses only deletion of
labels and right translations, so it applies to any quandle. This is the half-product
transported by $\theta^+$. For example, for degree-one cochains it gives
$(f\prec g)(x,y)=f(x)g(y)$.

For the remainder of this section, assume as in Section~\ref{sec:james-comparison} that $X$ and
$G$ are finite and $p\nmid |H|$. We take homology over a field $k$ of characteristic $p$ and
suppress $k$ from the notation. Put $U\coloneqq\widetilde H_*(V)$. Under the suspension
isomorphism $U\cong\susp H_*(\mathcal K)$ for $V=\Sigma\mathcal K_+$, write $\susp b\in U$ for
the suspension of $b\in H_*(\mathcal K)$.

Let $\pi:M\to BX$ be the covering quotient. Define the degree-$-1$ forgetful operator $\ell$ on
the shadow complex by
\[
 \ell(y;x_1,\ldots,x_n)\coloneqq(x_1;x_2,\ldots,x_n)\quad(n\ge1),
 \qquad \ell|_{S_0}\coloneqq 0\,.
\]
The map $\ell$ is the composite of the cellular chain map of $\mathcal N\to BX$ with the
rack--shadow correspondence $\rho$, so $\partial_S\ell=-\ell\partial_S$.

\begin{lemma}\label{half:shadow}
The map $v_*$ is an isomorphism of right $H_*(M)$-modules. Define the degree-one isomorphism
$\sigma_{\mathrm F}$ and the degree-$-1$ operator $Q_{\mathrm F}$ by
\[
 \sigma_{\mathrm F}\coloneqq \rho_*^{-1}v_*:
 H_n(M)\xrightarrow{\cong}H_{n+1}^R(X),\qquad
 Q_{\mathrm F}\coloneqq v_*^{-1}\ell_*v_*\,.
\]
Then
\begin{align}
 \pi_*^+&=\sigma_{\mathrm F}Q_{\mathrm F},
 \label{half:piQ}\\
 Q_{\mathrm F}(ab)&=Q_{\mathrm F}(a)b+\varepsilon(a)Q_{\mathrm F}(b),
 \label{half:Fox}\\
 \overline\Delta_{\prec}(\sigma_{\mathrm F}a)
 &=\sum_{(a)}\sigma_{\mathrm F}(a_{(1)})\otimes
     \sigma_{\mathrm F}(Q_{\mathrm F}a_{(2)})\,.
 \label{half:shadow-diagonal}
\end{align}
Here $\pi_*^+$ denotes the positive-degree part of $\pi_*$, set to zero in degree zero, and
$\overline\Delta_{\prec}$ is the half-diagonal dual to the canonical half-product, with both
factors in positive degrees. We use $\sum_{(a)}a_{(1)}\otimes a_{(2)}=\Delta a$ for the diagonal
on $M$.
\end{lemma}
\begin{proof}
The right $M$-action \cite[Proposition 3, (22)]{Clauwens} commutes with the left $H$-action, so
the quotient map $v$ is right $M$-equivariant. Equation~\eqref{good:cover-transfer} with
$\Lambda=k$ gives the isomorphism $v_*$. The label formulas give $\ell v_\#=\rho\pi_\#$ on
chains, with both sides zero in degree zero. This proves \eqref{half:piQ}, and $Q_{\mathrm F}$
vanishes in degree zero. If the first factor in the product has positive degree, $\ell$ takes
its first label as the region label. If the first factor is a vertex, $\ell$ takes the first
label of the second factor as the region label. This gives the Fox rule \eqref{half:Fox} for the
unsigned forgetful operator \cite[Proposition 4, Remark 5, (17)]{Clauwens}.

Comparing Clauwens's splitting of the product by the first coordinate \cite[Proposition 11,
(48)]{Clauwens} with the half-diagonal of \cite[Proposition 26]{CFLM} using the unsigned pairing
gives the following identity for shadow cochains $f,g$.
\begin{equation}\label{half:cochain}
 (\rho^*f)\prec(\rho^*g)=\rho^*(f\smile\ell^*g)\,.
\end{equation}
Deleting the first coordinate from the left subset gives the subset for the shadow cup product,
and the action on labels agrees because $\ell^*g$ forgets the region label. The first coordinate
contributes no inversions, so the shuffle signs also agree. Dualizing \eqref{half:cochain} and
using naturality of the diagonal under $v$ and the K\"unneth isomorphism gives
\eqref{half:shadow-diagonal}. In the graded tensor expression
$(\sigma_{\mathrm F}\otimes\sigma_{\mathrm F})(\id\otimes Q_{\mathrm F})\Delta$, the two Koszul
signs from $|Q_{\mathrm F}|=-1$ and $|\sigma_{\mathrm F}|=1$ cancel.
\end{proof}

\subsection{The half-shuffle and operations on words}
The conjugates $k^{-1}tk$ in \eqref{half:kernel} give the letters used in the comparison.
Since $\mathcal K\hookrightarrow\mathcal G\to\mathcal K$ is the identity, we
regard $H_*(\mathcal K)$ as a subalgebra of $H_*(\mathcal G)$. Let $\omega$ be a continuous map
representing $M\xrightarrow{c_M}\Omega Y\xrightarrow{\Omega\chi}\Omega W\simeq\mathcal G$, and
write $\omega_*$ for the induced Hopf algebra isomorphism. We retain the notation $\sigma_{\mathrm F}$ for
$\sigma_{\mathrm F}\omega_*^{-1}$ and $Q_{\mathrm F}$ for $\omega_*Q_{\mathrm F}\omega_*^{-1}$.

The orbit projection $r:B(*,G,M)\to G\backslash M=BX$ is a homotopy equivalence. In the model $EG\times_G M$ it is the bundle with contractible fiber $EG$ associated to
the covering $\pi:M\to BX$. Write $i_M:M\to B(*,G,M)$ for the inclusion of the fiber
of $B(*,G,M)\to BG$ over its vertex.
Since $ri_M=\pi$, a homotopy inverse of $r$ carries $\pi$ to $i_M$ up to homotopy.
On this fiber the prism of Proposition~\ref{james:space} assigns the reverse of the bar edge
$\ell_a$ to $a\in M$. Hence $\beta\pi\simeq\delta\circ\operatorname{inv}_{\mathcal G}\circ\omega$. The
reversal of the entire loop is the same as in \eqref{half:alpha}, and on homology we obtain
\begin{equation}\label{half:bar}
 \beta_*\pi_*=\delta_*\operatorname{inv}_{\mathcal G*}\omega_*\,.
\end{equation}

\begin{lemma}\label{half:letters}
Let $u\in H_1(\mathcal G)$ be the image of the positively oriented fundamental class $[S^1]$
under the inclusion of generators $S^1\hookrightarrow P\hookrightarrow\mathcal G$, and put
$c_k\coloneqq Q_{\mathrm F}u\in H_0(\mathcal G)=k$. Restrict $\Gamma$ along the canonical
inclusion $V\hookrightarrow J(V)$ and put $z_b\coloneqq(\Gamma|_V)_*(\susp b)$ for
$b\in H_d(\mathcal K)$, $d\ge0$. Writing $S$ for the antipode of $H_*(\mathcal K)$, we have
\begin{gather}
 Q_{\mathrm F}\bigl(H_*(\mathcal K)\bigr)=0,
 \label{half:QB}\\
 z_b=\sum_{(b)}(-1)^{|b_{(1)}|}S(b_{(1)})u b_{(2)},\qquad
 \Delta z_b=z_b\otimes1+1\otimes z_b,\qquad
 Q_{\mathrm F}z_b=c_kb\,.
 \label{half:z}
\end{gather}
\end{lemma}
\begin{proof}
The restriction $\delta|_{\mathcal{K}}$ is nullhomotopic by the definition of the path model of
the homotopy fiber, and inversion preserves $\mathcal K$. Thus \eqref{half:piQ} and
\eqref{half:bar} imply \eqref{half:QB}. Express right conjugation as the diagonal, inversion,
and a threefold product, and orient $\susp b$ by $[S^1]\times b$. Moving the first $\mathcal K$
factor past $S^1$ gives the sign $(-1)^{|b_{(1)}|}$, and inversion gives $S(b_{(1)})$. The class
$z_b$ is primitive because it comes from the reduced homology of the suspension $V$. Finally,
\eqref{half:Fox} and \eqref{half:QB} give
$Q_{\mathrm F}z_b=c_k\sum_{(b)}\varepsilon(b_{(1)})b_{(2)}=c_kb$.
\end{proof}

Concatenation of James words and the K\"unneth formula over a field give $H_*(J(V))=T(U)$
\cite[(1.12), \S3]{James}. Define the half-diagonal that places the first letter on the left by
\begin{equation}\label{half:unshuffle}
 \Delta^J_{\prec}(v_1\cdots v_m)
 \coloneqq \sum_{\substack{I\sqcup J=\{2,\ldots,m\}\\J\ne\varnothing}}
 \epsilon(I,J)(v_1v_I)\otimes v_J\,.
\end{equation}
Each subword is read in its original order, and $\epsilon(I,J)$ is the Koszul sign of the
separation into left and right subwords, given by
$\epsilon(I,J)=(-1)^{\sum_{j\in J,\,i\in I,\,j<i}|v_j||v_i|}$. When
$\Delta^J_{\prec}w=\sum w'\otimes w''$, define $(f\prec_J g)(w)\coloneqq \sum f(w')g(w'')$ using
the same unsigned pairing as for rack cochains. Let $U^\vee$ denote the degreewise dual of $U$. A
graded Zinbiel algebra is a graded module with a degree-zero bilinear operation $\prec$
satisfying
\[
(a\prec b)\prec c=a\prec(b\prec c+(-1)^{|b||c|}c\prec b)
\]
for homogeneous elements \cite[\S1.2]{Loday}.

\begin{theorem}[The rack half-product is the half-shuffle]\label{half:main}
Suppose $X$ and $G$ are finite and $p\nmid |H|$. For every field $k$ of characteristic $p$, the
isomorphism $\mathcal T_k$ in \eqref{half:coefficient-comparison} satisfies, for positive-degree
classes $a,b$,
\begin{equation}\label{half:equality}
 \mathcal T_k(a\prec b)=\mathcal T_k(a)\prec_J\mathcal T_k(b)\,.
\end{equation}
Consequently $H_R^{>0}(X;k)$ is isomorphic to the free graded Zinbiel algebra on $U^\vee$.
\end{theorem}
\begin{proof}
Put $\mathcal F_k\coloneqq \beta_*^{-1}\alpha_*$ and let $v_i\coloneqq \susp b_i$ for
homogeneous $b_i$. Write $Z_I$ for the ordered product of the $z_{b_i}$ with indices in $I$,
with $Z_\varnothing\coloneqq 1$. Since \eqref{half:alpha} and \eqref{half:bar} use the same
reversal of the entire loop and $\Gamma$ is multiplicative, \eqref{half:piQ}, \eqref{half:Fox},
and \eqref{half:z} give, for $m\ge1$,
\begin{equation}\label{half:word}
 \mathcal F_k(v_1\cdots v_m)=\pi_*\omega_*^{-1}Z_{\{1,\ldots,m\}}
 =\sigma_{\mathrm F}Q_{\mathrm F}Z_{\{1,\ldots,m\}}
 =c_k\sigma_{\mathrm F}(b_1Z_{\{2,\ldots,m\}})\,.
\end{equation}
Since the $z_{b_i}$ are primitive, their coproducts expand as sums over complementary
subwords. Substituting in \eqref{half:shadow-diagonal} gives
\[
 \overline\Delta_{\prec}\mathcal F_k(v_1\cdots v_m)
 =\sum_{\substack{I\sqcup J=\{2,\ldots,m\}\\J\ne\varnothing}}
 \epsilon(I,J)\,\mathcal F_k(v_1v_I)\otimes\mathcal F_k(v_J)\,.
\]
Indeed, the term $Q_{\mathrm F}(b_{1(2)}Z_J)$ in the right factor vanishes for $J=\varnothing$
by \eqref{half:QB}. If $j=\min J$ exists, the Fox rule gives
$\eps(b_{1(2)})c_k b_j Z_{J\setminus\{j\}}$. For terms that survive the counit, $|b_{1(2)}|=0$,
so passing this factor across $Z_I$ contributes $(-1)^{|b_{1(2)}||Z_I|}=1$. Only the sign in
\eqref{half:unshuffle} remains. The counit identity leaves $b_1Z_I$ on the left. Absorbing the
two factors $c_k$ into the two subword expressions in \eqref{half:word} gives the displayed
identity. For $I=\varnothing$ the left subword consists of $v_1$. For $m=1$ both sides vanish.

Since $X$ is finite and $\mathcal F_k$ is an isomorphism, $T(U)$ is finite dimensional in each
degree. As $U$ is concentrated in positive degrees, taking degreewise duals gives
\eqref{half:equality}. The free Zinbiel algebra description of \cite[Proposition 1.8]{Loday} applies with the
Koszul sign rule for permutations in the graded setting.
\end{proof}

Under the James stable splitting \cite{James}~\cite[Propositions 4.18, 5.8]{DevalapurkarHaine}
\begin{equation}\label{loop:stable-words}
 \Sigma J(V)\simeq\bigvee_{r\ge1}\Sigma V^{\wedge r}
\end{equation}
the suspension of the canonical inclusion $V\hookrightarrow J(V)$ corresponds to the inclusion of the word-length-one summand. Stable operations preserve each word-length summand. For $k=\F_p$, their values on words in the tensor description follow by naturality with respect to concatenation $V^r\to J(V)$, using the Cartan formulas and sign rules for external products for the Steenrod operations and the first Bockstein $\beta_1$ (Section~\ref{sec:operations}). Their values on letters follow from the suspension isomorphism for $V=\Sigma \mathcal{K}_+$.
\Needspace{7\baselineskip}

\section{Divided powers and Pontryagin operations}\label{sec:operations}
This section tests
commutation with the ordinary IK map on liftable classes of $A_X$, and the realization
after suspension gives commutation with stable operations in positive degrees.

\subsection{The canonical half-product and divided powers}

Let $B$ be a unital commutative coefficient ring. The construction of $\prec_Q$ and
Proposition~\ref{prop:divided-powers} apply to any quandle.
Corollary~\ref{half:divided-powers} uses the hypotheses of Theorem~\ref{half:main}. The half-diagonals and their chain
homotopies \cite[Propositions 26--27, Theorem 28]{CFLM} make positive-degree rack cohomology a
graded Zinbiel algebra under the canonical half-product $\prec$ of
Section~\ref{sec:loop-structures}. The other half-product, which sends the edge in the first
coordinate to the right factor, gives $(-1)^{|u||v|}v\prec u$ on cohomology. The sum of the two
half-products gives the ordinary cup product.
\begin{equation}\label{eq:half-cup-symmetrization}
 u\smile v=u\prec v+(-1)^{|u||v|}v\prec u\,.
\end{equation}
For any quandle $X$, the integral chain splitting of \cite[Proposition 30, Lemma 33]{CFLM}
gives $H_R^*(X;B)=\pi_X^*H_Q^*(X;B)\oplus H_D^*(X;B)$, naturally in $X$ and $B$.
In positive degrees the degenerate summand is a Zinbiel ideal
\cite[Propositions 27, 30, 32 and Theorem 34]{CFLM}. Indeed, both integral half-diagonals
send the image of the chain section into the tensor product of the rack complex with that
image. A $B$-valued cocycle in the degenerate summand can be chosen to vanish on this
image. With that cocycle in the right slot, either half-product again vanishes there.
The signed interchange of the two half-products on cohomology gives the assertion for
the left slot. Evaluation uses multiplication in $B$, and the degreewise sign
$\theta^+$ of Section~\ref{sec:loop-structures} converts to our cochain convention. Via the splitting, the quotient identifies with $H_Q^{>0}(X;B)$. Writing brackets
for projection to this quotient, set
\begin{equation}\label{eq:quandle-half-product}
 a\prec_Q b\coloneqq [\pi_X^*a\prec\pi_X^*b]\,.
\end{equation}
This quotient makes $H_Q^{>0}(X;B)$ a graded Zinbiel algebra. Since $\pi_X^*$ preserves the
ordinary cup product, the symmetrization of $\prec_Q$ is the ordinary cup product on $B^QX$.

For example, take the integral degree-one class $f(x)=1$ of the one-element quandle $X=\{x\}$.
Since $H_Q^2(X;\Z)=0$, we have $f\prec_Q f=0$. On the rack side the differential vanishes and
$(\pi_X^*f\prec\pi_X^*f)(x,x)=1$, so $\pi_X^*f\prec\pi_X^*f$ generates $H_R^2(X;\Z)\cong\Z$.

\begin{proposition}[Divided powers in even degrees]\label{prop:divided-powers}
For any unital commutative ring $B$, the positive even-degree parts of rack and quandle
cohomology are nonunital commutative algebras with divided powers under the ordinary cup
product. The divided powers are defined by $\gamma_1(x)\coloneqq x$ and
$\gamma_{r+1}(x)\coloneqq x\prec\gamma_r(x)$, using $\prec_Q$ in quandle cohomology.
Homomorphisms of coefficient rings and of quandles preserve these structures. In particular,
for every $r\ge1$,
\begin{equation}\label{eq:zinbiel-factorial}
 x^r=r!\,\gamma_r(x)\,.
\end{equation}
\end{proposition}
\begin{proof}
The sign in \eqref{eq:half-cup-symmetrization} is $+1$ in even degrees, so the theorem of Dokas
\cite[Proposition 3.3, Theorem 3.4]{Dokas} applies. Naturality follows from the definition
by iteration of the half-product.
\end{proof}

For a prime $p$, every positive even-degree rack or quandle class $x$ over $\F_p$ satisfies
$x^p=0$ by \eqref{eq:zinbiel-factorial}. For odd $p$ and $|x|=2m$, the instability relation
gives $P^m(x)=x^p=0$. Over $\F_2$, every homogeneous positive-degree rack or quandle class $x$
satisfies $\operatorname{Sq}^{|x|}x=x^2$, and \eqref{eq:half-cup-symmetrization} and its quandle
analogue give $x^2=0$.

\begin{corollary}[Word formulas over a field]\label{half:divided-powers}
Under the assumptions of Theorem~\ref{half:main}, for every positive even-degree rack class $x$
and $m\ge1$ we have
\[
 \mathcal T_k\gamma_m(x)=\gamma_m^J(\mathcal T_kx)\,.
\]
Here $\gamma_m^J$ is the divided power defined using $\prec_J$. Taking the degreewise dual of
the words in $H_*(J(V);k)\cong T(U)$ identifies the target of $\mathcal T_k$ with $\operatorname{Sh}(U^\vee)$, the shuffle algebra on $U^\vee$ with Koszul signs. There is then an
isomorphism of graded Zinbiel algebras
\[
 (H_Q^{>0}(X;k),\prec_Q)
 \cong(\operatorname{Sh}^{>0}(U^\vee),\prec_J)/\mathcal T_k(H_D^{>0}(X;k))\,.
\]
\end{corollary}
\begin{proof}
Apply Theorem~\ref{half:main} to the iterated half-products and pass to the quotient by
$H_D^{>0}$. The choice of a word basis for the quandle quotient depends on the explicit form of
this ideal.
\end{proof}

\subsection{The IK map and cohomology operations}

\begin{proposition}[Realization after suspension]\label{ik:stable}
Let $X$ be a connected discrete quandle with a transitive augmentation whose image generates
$G$. There is a continuous map $\Sigma BX\to\Sigma A_X$ which, under the suspension
isomorphisms, induces $\pi_X^*\Phi_B$ in positive degrees for every constant coefficient module
$B$. Thus the map to rack cohomology commutes in positive degrees with every ordinary stable
cohomology operation $\mathcal O$ between constant coefficient modules $B,C$. The corresponding
identity in quandle cohomology is
\begin{equation*}
 \mathcal O\Phi_B=\Phi_C\mathcal O\,.
\end{equation*}
In particular, Steenrod operations and Bocksteins commute with the ordinary IK map.
\end{proposition}
\begin{proof}
Use the map $S^1\times\mathcal N\to BX$ of Fenn, Rourke, and Sanderson~\cite[Theorem 5.14
and the paragraph following it]{FRS}. With the circle factor first, the product of its
positively oriented $1$-cell with the shadow cell $(y;\boldsymbol x)$ maps to the rack cell
$(y,\boldsymbol x)$ with its given orientation. Restrict a homotopy inverse of the standard
splitting $\Sigma(S^1_+)\simeq S^1\vee S^2$ to the $S^2$ summand and choose a cellular
representative $S^2\to\Sigma(S^1_+)$ sending the positive fundamental class to the suspension of
the fundamental class of the circle. Smash this map with the identity of $\mathcal N_+$. Using
$S^1_+\wedge\mathcal N_+\cong(S^1\times\mathcal N)_+$, follow it with the suspension of the
circle map and the basepoint quotient $\Sigma BX_+\to\Sigma BX$ to obtain
$\upsilon:\Sigma^2\mathcal N_+\to\Sigma BX$. Thus $\upsilon_\#(\susp^2c)=\susp\rho^{-1}c$ on
integral cellular chains. Hence $\upsilon$ induces an integral homology isomorphism and is a
homotopy equivalence, since the spaces are simply connected and have CW type.

Take the auxiliary vertex $p=x_0$. The quotient $\mathcal N_+\to\mathcal N$ identifies the
disjoint basepoint with $x_0$ and
sends the reduced $0$-chain $y-*$ to $y-x_0$. We have $\Theta_0(x_0)=q(x_0,x_0)=0$ by simplicial
degeneracy, and $\varphi=\Theta\rho$ in positive degrees by \eqref{eq:ordinaryIK}.
Theorem~\ref{thm:comparison} therefore shows that the following composite, using a homotopy
inverse of $\upsilon$, represents $\susp\varphi$ up to integral chain homotopy.
\[
 \Sigma BX\xrightarrow{\ \upsilon^{-1}\ }\Sigma^2\mathcal N_+
 \longrightarrow\Sigma^2\mathcal N\xrightarrow{\ \Sigma\widehat\kappa\ }\Sigma A_X\,.
\]
On rack cohomology the identity follows from naturality and compatibility of stable operations
with suspension. On quandle cohomology it follows from naturality under $\pi_X$ and injectivity
of $\pi_X^*$ for the output coefficients $C$.
\end{proof}

We call the generalized Pontryagin operations of Browder and Thomas~\cite[Theorems
1--2]{BrowderThomas} Pontryagin--Thomas operations. For a prime $p$ and integers $r,m\ge1$,
they give an unstable primary cohomology operation natural in the space $Z$,
\[
 C_{p,r,Z}:H^{2m}(Z;\Z/p^r)\longrightarrow
 H^{2mp}(Z;\Z/p^{r+1})\,.
\]
For $p=2$, Theorem~2 of the same paper gives the generalized Pontryagin square. For odd $p$ and
$r=1$,
it is the operation that Clauwens calls the Thomas operation \cite[\S6.3, Proposition
43]{Clauwens}. When $Z=B^QX$, we write $C_{p,r,Q}$. Let $\eta_r:\Z/p^{r+1}\to\Z/p^r$ be
coefficient reduction and let $\iota_r:\Z/p^r\to\Z/p^{r+1}$ be given by $[a]\mapsto[pa]$. We use
the same notation for the induced maps on cohomology. The lifting formula of Browder and Thomas~\cite[(1.3)]{BrowderThomas} is
\begin{equation}\label{eq:pontryagin-lift}
 C_{p,r,Z}(\eta_r\widetilde x)=\widetilde x^p,\qquad
 \widetilde x\in H^{2m}(Z;\Z/p^{r+1})\,.
\end{equation}

We use the following two short exact sequences of coefficients.
\begin{align*}
 0&\longrightarrow\F_p\xrightarrow{\ p^r\ }\Z/p^{r+1}
   \xrightarrow{\ \eta_r\ }\Z/p^r\longrightarrow0\,,\\
 0&\longrightarrow\Z/p^r\xrightarrow{\ \iota_r\ }\Z/p^{r+1}
   \longrightarrow\F_p\longrightarrow0\,.
\end{align*}
A class lifts modulo $p^{r+1}$ exactly when it lies in the kernel of the connecting
homomorphism for the first sequence. For the second sequence, write $\delta_r$ for the
connecting homomorphism. Then $\ker\iota_r=\im\delta_r$ on cohomology, which tests whether a
divided power maps to zero under $\iota_r$. For $r=1$ both sequences give the first Bockstein
$\beta_1$.

\begin{proposition}[Divided powers and the IK map in the coefficient tower]\label{prop:coefficient-tower}
Let $X$ be any quandle. If $x\in H_Q^{2m}(X;\Z/p^r)$ lifts to $\Z/p^{r+1}$ coefficients, then
\begin{equation}\label{eq:tower-halfpower}
 C_{p,r,Q}(x)=(p-1)!\,\iota_r\gamma_p(x)\,.
\end{equation}
The nonvanishing criterion is
\begin{equation}\label{eq:tower-nonzero}
 C_{p,r,Q}(x)\ne0\quad\Longleftrightarrow\quad
 \gamma_p(x)\notin\im\bigl(\delta_r:H_Q^{2mp-1}(X;\F_p)
 \longrightarrow H_Q^{2mp}(X;\Z/p^r)\bigr)\,.
\end{equation}
The same formula and criterion hold in rack cohomology.

Suppose further that $X$ satisfies the assumptions of Section~\ref{sec:ik-comparison}. If
$a\in H^{2m}(A_X;\Z/p^r)$ lifts to $\Z/p^{r+1}$ coefficients and $x=\Phi_{\Z/p^r}(a)$, then
\begin{equation}\label{eq:tower-ik-defect}
 \Phi_{\Z/p^{r+1}}\bigl(C_{p,r,A_X}(a)\bigr)=0\,.
\end{equation}
The class $C_{p,r,Q}(x)$ is given by \eqref{eq:tower-halfpower}. Thus whenever the condition in
\eqref{eq:tower-nonzero} holds, the original IK map fails to commute with this generalized
Pontryagin operation.
\end{proposition}
\begin{proof}
Choose $\widetilde x$ with $\eta_r\widetilde x=x$. By \eqref{eq:pontryagin-lift},
\eqref{eq:zinbiel-factorial}, and $\iota_r\eta_r=p\,\id$,
\[
 C_{p,r,Q}(x)=\widetilde x^p=p!\gamma_p(\widetilde x)
 =(p-1)!\,\iota_r\gamma_p(x)\,.
\]
The last equality uses naturality of divided powers under $\eta_r$. Since $(p-1)!$ is a unit
modulo $p^{r+1}$,
the equality $\ker\iota_r=\im\delta_r$ gives the criterion for nonvanishing. The rack case is
identical. The IK image of a lift $\widetilde a$ of $a$ lifts $x$. The IK image of
$C_{p,r,A_X}(a)=\widetilde a^p$ is zero by Corollary~\ref{lem:cupzero}.
\end{proof}

\begin{corollary}[The word formula for mod-$p$ input]\label{prop:pontryagin-halfpower}\label{half:pontryagin-first}
Let $X$ be any quandle and let $x\in H_Q^{2m}(X;\F_p)$ with $m\ge1$. If $\beta_1x=0$, then
\begin{equation*}
 C_{p,1,Q}(x)=-\iota_1\gamma_p(x)
 \quad\text{in }H_Q^{2mp}(X;\Z/p^2)\,.
\end{equation*}
For $p=2$ the sign on the right can be dropped because $2\iota_1=0$. The nonvanishing criterion
is \eqref{eq:tower-nonzero} with $r=1$. Under the assumptions of Theorem~\ref{half:main}, take the James comparisons for $\F_p$
and $\Z/p^2$ to be induced by the same maps \eqref{james:zigzag}. Then
\[
 \mathcal T_{\Z/p^2}\bigl(\pi_X^*C_{p,1,Q}(x)\bigr)
 =-\iota_1\gamma_p^J\bigl(\mathcal T_{\F_p}(\pi_X^*x)\bigr)\,.
\]
\end{corollary}
\begin{proof}
The first formula follows from $\ker\beta_1=\im\eta_1$,
Proposition~\ref{prop:coefficient-tower}, and Wilson's theorem $(p-1)!\equiv-1\pmod p$. For the
James formula, $\pi_X^*x$ also lifts to $\Z/p^2$ coefficients because $\pi_X^*$ commutes with
coefficient reduction. Naturality gives $\pi_X^*C_{p,1,Q}(x)=C_{p,1,BX}(\pi_X^*x)$. Apply the rack version of Proposition~\ref{prop:coefficient-tower} to the right-hand side.
Corollary~\ref{half:divided-powers} gives the word formula for the divided power over
$\F_p$. The James comparison commutes with $\iota_1$, which carries that formula to the
output with $\Z/p^2$ coefficients.
\end{proof}
\Needspace{7\baselineskip}

\section{Cup powers and nonvanishing of Pontryagin operations}\label{james:section}
In this section, we study the cup powers and Pontryagin operations in more details.

\subsection{Primitive extensions and low-degree homology}

Let $p$ be a prime, put $R\coloneqq \Z_{(p)}$, and write $\nu_p$ for the $p$-adic valuation. Let
$C_m$ denote the cyclic group of order $m$. We write $\bmod p^s$ for coefficient reduction from
$R$ to $\Z/p^s$ and take the additive order of zero to be $1$. For a space $Z$ and $s\ge1$, the
coefficient sequence $0\to R\xrightarrow{p^s}R\to\Z/p^s\to0$ shows that the kernel of reduction
on $H^n(Z;R)$ is $p^sH^n(Z;R)$. Thus the order of a reduced class is its order in the quotient
$H^n(Z;R)/p^sH^n(Z;R)$.

\begin{lemma}[Evaluation with prime-power coefficients]\label{james:mod4}
For every space $Z$ and integer $s\ge1$, the natural evaluation map
\[
 H^n(Z;\mathbb Z/p^s)\longrightarrow
 \operatorname{Hom}_{\mathbb Z/p^s}\bigl(H_n(Z;\mathbb Z/p^s),\mathbb Z/p^s\bigr)
\]
is an isomorphism. Consequently any class of order $p^s$ on the left evaluates to $1$ on some homology class.
\end{lemma}
\begin{proof}
Put $S=\Z/p^s$. The map $1\mapsto 1/p^s+\Z$ identifies its additive group with the subgroup of $\mathbb Q/\Z$ killed by $p^s$.
For every $S$-module $M$, this identifies $\Hom_S(M,S)$ with
$\Hom_{\Z}(M,\mathbb Q/\Z)$. The latter functor is exact since $\mathbb Q/\Z$ is divisible
and hence injective. Dualizing $C_*(Z;S)$ gives the asserted evaluation isomorphism.
A functional of order $p^s$ has image $S$ and hence takes the value $1$.
\end{proof}

\begin{lemma}[Primitive James extensions and evaluation of powers]\label{james:quadratic}
Let $V$ be a connected CW complex whose basepoint is a $0$-cell, let $n>0$ be even, and let
$j:V\to J(V)$ be the canonical inclusion. Every $a\in\widetilde H^n(V;R)$ extends to a primitive
class $\widehat a\in H^n(J(V);R)$ with $j^*\widehat a=a$. Here primitive means that the pullback of $\widehat a$ under James multiplication equals
the sum of its pullbacks under the two projections. If $a$ has
finite order, $\widehat a$ may be chosen to have the same order. For $s,m\ge1$, put
$a_s\coloneqq a\bmod p^s$ and $\widehat a_s\coloneqq \widehat a\bmod p^s$. For the $m$-fold
James multiplication $\mu_m$ and $c\in H_n(V;\Z/p^s)$, we have
\begin{equation}\label{james:factorial-evaluation}
 \bigl\langle\widehat a_s^m,\mu_{m*}((j_*c)^{\times m})\bigr\rangle
 =m!\,\langle a_s,c\rangle^m\quad\text{in }\Z/p^s\,.
\end{equation}
\end{lemma}
\begin{proof}
Take $K(R,n)$ to be the topological abelian group obtained from $R$ in degree $n$ by the
Dold--Kan correspondence and realization \cite[Ch.~III, Corollary 2.3]{GoerssJardine}. By the defining property of $J(V)$, a based
representative $f:V\to K(R,n)$ of $a$ extends to the monoid map
$\widehat f(v_1\cdots v_r)=\sum_i f(v_i)$. This construction respects based homotopies and
addition, so $a\mapsto\widehat a=[\widehat f]$ is an additive section of $j^*$ and preserves
orders. The word formula gives $\mu_m^*\widehat a=\sum_i\operatorname{pr}_i^*\widehat a$.
After reduction modulo $p^s$, only terms using each factor once can contribute to the pairing
of $\mu_m^*(\widehat a_s^m)$ with $(j_*c)^{\times m}$. Each of the $m!$ such terms contributes
$\langle a_s,c\rangle^m$ with positive sign since $n$ is even. This proves
\eqref{james:factorial-evaluation}.
\end{proof}

We apply this evaluation under a low-degree vanishing hypothesis. Let $X$ be a finite
connected quandle with a transitive augmentation whose image generates a finite group $G$,
and let $H$ be a point stabilizer. Assume $p\nmid |H|$, fix an even integer $N\ge4$, and
suppose that
\begin{equation}\label{james:criterion-hypotheses}
 H_j(A_X;R)=0\qquad(2\le j\le N-2)\,.
\end{equation}
Use $\mathcal{K}=\Omega A_{X,R}$ and $V=\Sigma \mathcal{K}_+$ from
Proposition~\ref{james:space}, put $T\coloneqq J(V)$, and let $j:V\to T$ be the canonical
inclusion.

\begin{lemma}[The IK isomorphism and primitive classes in degree $N$]\label{james:degree4}
There are isomorphisms
\[
 \Phi_R^N:H^N(A_X;R)\xrightarrow{\ \cong\ }H_Q^N(X;R),\qquad
 H^N(A_X;R)\cong\Ext_R^1(H_{N-1}(A_X;R),R)\,.
\]
In degree $N$, the restriction to the subspace $V$ of words of length one determines each
torsion class in the James model. Both maps below are isomorphisms, and every class in
$H^N(T;R)_{\mathrm{tors}}$ is primitive for James multiplication.
\begin{equation}\label{james:degree4-comparison}
 H_Q^N(X;R)\xrightarrow{\ \mathcal T_R\pi_X^*\ }H^N(T;R)_{\mathrm{tors}}
 \xrightarrow{\ j^*\ }H^N(V;R)\,.
\end{equation}
Let $s\ge1$ and suppose that the reduction $a\bmod p^s$ of $a\in H^N(A_X;R)$ has order $p^s$. If
$y_s\coloneqq\Phi_R^N(a)\bmod p^s$, then for every $m\ge1$,
\[
 \operatorname{ord}(y_s^m)=p^{\max(s-\nu_p(m!),0)}\,.
\]
\end{lemma}
\begin{proof}
Put $B\coloneqq H_{N-1}(A_X;R)$. The map $q_G$ of Lemma~\ref{good:marked-pinch}, the homology of
finite groups, and the rational homology of finite connected quandles \cite{EtingofGrana,LN}
show that $H_i(A_X;R)$ for $i>0$ and $H_i^Q(X;R)$ for $i\ge2$ are finite $p$-groups. Both chain
complexes are finitely generated in each degree. The $\Hom_R$ terms in the universal coefficient
theorem in degree $N$ vanish on both sides. By naturality, $\Phi_R^N$ agrees with the
isomorphism induced by applying $\Ext_R^1(-,R)$ to the degree-$(N-1)$ homology isomorphism in
Theorem~\ref{good:degree3}. Thus $H^N(A_X;R)$ and $B=H_{N-1}(A_X;R)$ have the same cyclic
invariant factors, since $\Ext_R^1(R/p^d,R)\cong R/p^d$.

By the proof of Theorem~\ref{good:degree3} and \eqref{james:criterion-hypotheses}, the space
$A_{X,R}$ is $(N-2)$-connected and $\mathcal K$ is $(N-3)$-connected. In the homology Serre
spectral sequence of the path fibration, $E^2_{i,j}=0$ for $i,j>0$ and $i+j<2N-3$. Since the
total space is contractible, the transgressions $d_{N-1}$ and $d_N$ from base degrees $N-1,N$ to
fiber degrees $N-2,N-1$, respectively, are isomorphisms. Thus
\begin{equation}\label{james:criterion-loop-low}
 H_i(\mathcal{K};R)=0\ (1\le i<N-2),\quad H_{N-2}(\mathcal{K};R)\cong B,\quad H_{N-1}(\mathcal{K};R)\cong H_N(A_X;R)\,.
\end{equation}
In degrees at most $N-1$, the reduced $R$-homology of $V=\Sigma\mathcal{K}_+$ is $R$ in
degree $1$, $B$ in degree $N-1$, and zero in all other degrees. Moreover $H_N(V;R)$ is
finite. The stable splitting
\eqref{loop:stable-words} and the K\"unneth formula give $H_{N-1}(T;R)=R\oplus B$. The free
summand comes from the word of length $N-1$ in the degree-one letter, and the $B$ summand is the
image of the word-length-one inclusion $j$. By naturality of the universal coefficient theorem,
the restriction of $j^*$ to torsion is the map $\Ext_R^1(R\oplus B,R)\to\Ext_R^1(B,R)$ induced
by this inclusion. Thus $j^*:H^N(T;R)_{\mathrm{tors}}\to H^N(V;R)\cong\Ext_R^1(B,R)$ is an
isomorphism. The first map in \eqref{james:degree4-comparison} is injective by the injectivity
of $\pi_X^*$ and Proposition~\ref{james:space}. Its source and target have the same order, so it
is an isomorphism. The primitive extension of Lemma~\ref{james:quadratic} preserves torsion and is the inverse
of this restriction map.

Put $\theta\coloneqq \mathcal T_R\pi_X^*\Phi_R^N(a)$. The map $a\mapsto j^*\theta$ is an
isomorphism, and on both $A_X$ and $V$ the kernel of reduction is the subgroup of
$p^s$-multiples. Hence $j^*\theta\bmod p^s$ has order $p^s$. By Lemma~\ref{james:mod4}, choose
$c$ with $\langle j^*\theta\bmod p^s,c\rangle=1$. Equation~\eqref{james:factorial-evaluation}
then shows that $(\theta\bmod p^s)^m$ evaluates to $m!$ on $\mu_{m*}((j_*c)^{\times m})$. The
comparison in Proposition~\ref{james:space} and $\pi_X^*$ are compatible with coefficient change
and products, so the order of $y_s^m$ is at least the order of this evaluation. Since
$\gamma_m(y_s)$ has coefficients in $\Z/p^s$, the identity $y_s^m=m!\,\gamma_m(y_s)$ gives the
matching upper bound.
\end{proof}

\subsection{Orders of cup powers and Pontryagin--Thomas values}

We first give a sufficient condition for the operation on the mod-$p^r$ reduction of an IK image
over $R$ to have order $p^r$. We then fix a generator of a cyclic direct summand and compute the
exact order of each power of its IK image.

\begin{theorem}[A criterion in terms of low-degree homology]\label{james:ik-criterion}
Under the assumptions of Lemma~\ref{james:degree4}, let $r\ge1$ and suppose that the
reduction of $a\in H^N(A_X;R)$ modulo $p^{r+1}$ has additive order $p^{r+1}$. Put $y_R\coloneqq \Phi_R^N(a)$ and
$y_r\coloneqq y_R\bmod p^r$. Then
\begin{align*}
 C_{p,r,Q}(y_r)&=(p-1)!\,\iota_r\gamma_p(y_r)\ne0,
 \\
 \Phi_{\Z/p^{r+1}}^{Np}\bigl(C_{p,r,A_X}(a\bmod p^r)\bigr)&=0
 \,.
\end{align*}
The nonzero value $C_{p,r,Q}(y_r)$ has additive order exactly $p^r$. Such an $a$ exists if and
only if the finite $p$-group $H_{N-1}(A_X;R)$ has exponent at least $p^{r+1}$.
\end{theorem}
\begin{proof}
By \eqref{eq:pontryagin-lift}, $C_{p,r,Q}(y_r)=(y_R\bmod p^{r+1})^p$.
Lemma~\ref{james:degree4} with $s=r+1$ and $m=p$ gives its order as
$p^{r+1-\nu_p(p!)}=p^r$. Equations~\eqref{eq:tower-halfpower} and \eqref{eq:tower-ik-defect}
give the divided power formula and the vanishing of the IK image of $C_{p,r,A_X}(a\bmod p^r)$.
By Lemma~\ref{james:degree4}, $H^N(A_X;R)$ and $H_{N-1}(A_X;R)$ have the same cyclic
invariant factors. Since the kernel of reduction is $p^{r+1}H^N(A_X;R)$, the required $a$
exists exactly when one of these factors has order at least $p^{r+1}$.
\end{proof}

\begin{theorem}[Orders of cup powers of an IK image]\label{james:all-powers}
In the setting of Lemma~\ref{james:degree4}, suppose $a\in H^N(A_X;R)$ generates a cyclic direct
summand of order $p^e$, where $e\ge1$. Put $y_R\coloneqq \Phi_R^N(a)$. For all integers
$m,s\ge1$ and $b\ge0$,
\begin{align}
 \operatorname{ord}\bigl((p^b y_R)^m\bigr)&=p^{\max(e-\nu_p(m!)-bm,0)},
 \label{james:integral-power-order}\\
 \operatorname{ord}\bigl((p^b y_R)^m\bmod p^s\bigr)
 &=p^{\max(\min(s,e)-\nu_p(m!)-bm,0)}
 \label{james:coefficient-power-order}\,.
\end{align}
Thus, if $N_p(e)\coloneqq \min\{m\ge1:\nu_p(m!)\ge e\}$, the nilpotence index of $y_R$ is
exactly $N_p(e)$, and
\begin{equation}\label{james:generated-ring}
 R[y_R]\cong R[t]/\bigl(p^{\max(e-\nu_p(m!),0)}t^m:1\le m\le N_p(e)\bigr),
 \qquad |t|=N\,.
\end{equation}
\end{theorem}
\begin{proof}
For coefficients modulo $p^s$ set $d\coloneqq \min(s,e)$, and for coefficients in $R$ set
$d\coloneqq e$. Write $y$ for the image of $y_R$ in the given coefficients. Since $p^d y=0$,
bilinearity of the half-product gives $p^d\gamma_m(y)=0$. By \eqref{eq:zinbiel-factorial}, the
order of $y^m$ is at most $p^{\max(d-\nu_p(m!),0)}$. Since $a$ generates a cyclic direct
summand, its reduction modulo $p^d$ has order $p^d$. Lemma~\ref{james:degree4} gives the
matching lower bound already for the reduction of $y^m$ modulo $p^d$. The identity
$(p^by_R)^m=p^{bm}y_R^m$ gives the formula for general $b$. The powers belong to distinct
degrees, so their annihilator ideals together with $y_R^{N_p(e)}=0$ give
\eqref{james:generated-ring}.
\end{proof}

\begin{corollary}[Orders of Pontryagin--Thomas values]\label{james:entire-tower}\label{james:threshold}
Under the assumptions of Theorem~\ref{james:all-powers}, the additive order of $y_R^p$ is
exactly $p^{e-1}$. For each $r\ge1$, put $y_r\coloneqq y_R\bmod p^r$. Then
\begin{gather}
 C_{p,r,Q}(y_r)=(y_R^p)\bmod p^{r+1},\qquad
 \operatorname{ord}\bigl(C_{p,r,Q}(y_r)\bigr)=p^{\min(r,e-1)},
 \label{james:whole-tower-order}\\
 \Phi_{\Z/p^{r+1}}^{Np}\bigl(C_{p,r,A_X}(a\bmod p^r)\bigr)=0
 \notag\,.
\end{gather}
More generally, in the setting of Lemma~\ref{james:degree4}, there exists $a\in H^N(A_X;R)$ with
$\Phi_R^N(a)^p\ne0$ if and only if $H_{N-1}(A_X;R)$ has exponent at least $p^2$. If
$pH_{N-1}(A_X;R)=0$, then $C_{p,r,Q}(\Phi_R^N(a)\bmod p^r)=0$ for every $a\in H^N(A_X;R)$ and
$r\ge1$.
\end{corollary}
\begin{proof}
Apply \eqref{eq:pontryagin-lift} and \eqref{eq:tower-ik-defect} to reductions from $R$, and take
$m=p$ and $b=0$ in Theorem~\ref{james:all-powers}. Since $\nu_p(p!)=1$, this gives the stated
orders. If the exponent of $H_{N-1}(A_X;R)$ is at least $p^2$, the $\Ext_R^1$ description in
Lemma~\ref{james:degree4} gives a generator of a direct summand with $e\ge2$. If the exponent is
at most $p$, then $py_R=0$ for every degree-$N$ class $y_R=\Phi_R^N(a)$. Thus $p\gamma_p(y_R)=0$
by bilinearity of the half-product. It follows that $y_R^p=p!\gamma_p(y_R)=0$, so its reduction
at every coefficient step also vanishes.
\end{proof}

\begin{corollary}\label{james:nonrealization}
In the setting of Lemma~\ref{james:degree4}, suppose $H_{N-1}(A_X;R)$ has exponent at least
$p^2$. Then $\Phi_R$ fails to preserve cup products in positive degrees. Neither $BX$ nor $B^QX$
admits a continuous map to $A_X$ inducing $\Phi_R$ in degrees $N$ and $Np$, where the required
map on rack cohomology is $\pi_X^*\Phi_R$. For odd $p$, the weaker hypothesis $H_{N-1}(A_X;R)\ne0$
suffices. In this case $\Phi_R$ fails to preserve cup products, and the corresponding
realization is impossible already in degrees $N$ and $2N$.
\end{corollary}
\begin{proof}
Let $a\in H^N(A_X;R)$ generate a cyclic direct summand of order $p^e$. If $m\ge2$ and
$\nu_p(m!)<e$, Theorem~\ref{james:all-powers} gives $\Phi_R(a)^m\ne0$ and
Corollary~\ref{lem:cupzero} gives $\Phi_R(a^m)=0$. Pullback along a continuous map preserves cup
products, and $\pi_X^*$ is an injective ring homomorphism. This rules out the stated realization
in degrees $N,mN$ for either space. The $\Ext_R^1$ description in Lemma~\ref{james:degree4}
gives a cyclic direct summand with $e\ge2$ in the first case and $e\ge1$ in the second. Take
$m=p$ in the first case, where $\nu_p(p!)=1$, and $m=2$ in the second, where $\nu_p(2!)=0$
because $p$ is odd.
\end{proof}

\begin{corollary}[Cyclic direct summands in higher degrees]\label{james:high-torsion}
Under the assumptions of Theorem~\ref{james:all-powers}, let $m\ge1$ satisfy $\nu_p(m!)<e$. Then
$\gamma_m(y_R)$ generates a cyclic direct summand of order $p^e$ in $H_Q^{Nm}(X;R)$. Moreover,
$H_{Nm-1}^Q(X;R)$ also contains $C_{p^e}$ as a direct summand.
\end{corollary}
\begin{proof}
Put $\Gamma=\gamma_m(y_R)$, $\bar\Gamma=\Gamma\bmod p^e$ and $d=\nu_p(m!)<e$.
Since $p^e y_R=0$, the recursion for divided powers gives $p^e\Gamma=0$. Since $m!/p^d$ is a unit,
\eqref{eq:zinbiel-factorial} and \eqref{james:coefficient-power-order} give
$\operatorname{ord}(p^d\bar\Gamma)=p^{e-d}$. Thus $\bar\Gamma$ has order $p^e$.
By Lemma~\ref{james:mod4} applied to $B^QX$, some $c\in H_{Nm}^Q(X;\Z/p^e)$ satisfies
$\langle\bar\Gamma,c\rangle=1$. Hence $x\mapsto\langle x\bmod p^e,c\rangle$ is a left
inverse to $\Z/p^e\to H_Q^{Nm}(X;R)$, $1\mapsto\Gamma$. For $Nm\ge4$, the groups $H_{Nm}^Q(X;R)$ and $H_{Nm-1}^Q(X;R)$ are finite
$p$-groups, as in the first paragraph of the proof of Lemma~\ref{james:degree4}.
Thus the $\Hom_R$ term in the universal coefficient sequence vanishes, giving
$H_Q^{Nm}(X;R)\cong\Ext_R^1(H_{Nm-1}^Q(X;R),R)$. Since $\Ext_R^1(R/p^j,R)\cong R/p^j$,
the homology group has the same cyclic invariant factors and contains $C_{p^e}$ as a direct summand.
\end{proof}
\Needspace{7\baselineskip}

\section{Applications to finite connected quandles}\label{sec:binary}
We will show that, for each prime $p$, degree-four IK images over $\Z_{(p)}$ have unbounded nilpotence indices
in a family of finite symplectic quandles of fixed type. Cyclic Alexander quandles of type
$3$ give degree-six examples over $\Z_{(7)}$.

\subsection{Symplectic quandles over finite fields}

Let $\F_q$ be a finite field of characteristic $\ell$, with $q\ge3$. On the set of nonzero row
vectors define
\begin{equation}\label{symplectic:operation}
 X_q\coloneqq \F_q^2\setminus\{0\},\qquad v*w\coloneqq v+\det(v,w)w\,.
\end{equation}
This is the nonzero part of the two-dimensional symplectic quandle \cite[Definition
3]{NavasNelson}~\cite[Example 7.2]{NosakaHomotopical}.

\begin{proposition}\label{binary:quandle}
Put $G\coloneqq \mathrm{SL}_2(\F_q)$, $H\coloneqq \{E_{12}(t):t\in\F_q\}$, and
$h\coloneqq E_{12}(1)$, where
$E_{12}(t)\coloneqq \left(\begin{smallmatrix}1&t\\0&1\end{smallmatrix}\right)$. Under the
bijection $H\backslash G\to X_q$ sending $Hg$ to the second row of $g$, the quandle operation is
$Hg*Hk=Hgk^{-1}hk$. The image of this augmentation generates $G$. The quandle $X_q$ is
connected, has $q^2-1$ elements and type $\ell$, and has inner automorphism group $\Inn(X_q)=G$
with point stabilizer $H$.
\end{proposition}
\begin{proof}
The group $H$ is abelian and $h\in Z(H)$. The right action of $G$ on nonzero row vectors is
transitive, and $(0,1)$ has stabilizer $H$. If $(c,d)$ is the second row of $g$, then
$g^{-1}hg=I_2+(d,-c)^{\mathsf t}(c,d)$, so the right action gives \eqref{symplectic:operation}.

Conjugation by diagonal matrices shows that the normal subgroup generated by $h$ contains all
$E_{12}(t^2)$. For $\ell=2$, the squaring map is surjective. For $\ell\ne2$, the identity
$t=((t+1)/2)^2-((t-1)/2)^2$ shows that this subgroup again contains all $E_{12}(t)$. Conjugation
by a Weyl matrix gives the lower triangular elementary matrices as well, so the conjugates of
$h$ generate $G$. Since $X_q$ contains the standard basis, the action is faithful. Hence
$\Inn(X_q)=G$, and the type is the order $\ell$ of $h$.
\end{proof}

\begin{theorem}[Cup powers and Pontryagin--Thomas values for the symplectic quandle $X_q$]\label{binary:main}
Let $p$ be a prime dividing $q^2-1$, and put $e\coloneqq \nu_p(q^2-1)$ and
$R\coloneqq \Z_{(p)}$. Then
\begin{equation}\label{binary:target-homology}
 H_1(A_{X_q};R)=H_2(A_{X_q};R)=H_4(A_{X_q};R)=0,
 \qquad H_3(A_{X_q};R)\cong C_{p^e},
\end{equation}
and $\Phi_R^4:H^4(A_{X_q};R)\to H_Q^4(X_q;R)$ is an isomorphism between cyclic groups of order $p^e$. For a generator $a\in H^4(A_{X_q};R)$, put $y_R\coloneqq \Phi_R^4(a)$. For all $m,r\ge1$,
\begin{align}
 \operatorname{ord}(y_R^m)&=p^{\max(e-\nu_p(m!),0)},\label{symplectic:power-order}\\
 \operatorname{ord}C_{p,r,Q}(y_R\bmod p^r)&=p^{\min(r,e-1)},\qquad
 \Phi_{\Z/p^{r+1}}C_{p,r,A_{X_q}}(a\bmod p^r)=0\,.
 \label{symplectic:tower}
\end{align}
The orders after coefficient reduction and the presentation of the generated subring are given
by Theorem~\ref{james:all-powers}.
\end{theorem}
\begin{proof}
Since $p\nmid|H|=q$, the map $q_G$ of Lemma~\ref{good:marked-pinch} gives
$H_*(A_{X_q};R)\cong H_*(G;R)$. In positive degrees, $H_j(G;\Z[1/\ell])\cong C_{q^2-1}$ for
$j\equiv3\pmod4$ and vanishes in the other degrees \cite[Corollary 3.6]{Hutchinson}. Changing
coefficients to $R$ gives \eqref{binary:target-homology}. Lemma~\ref{james:degree4} with $N=4$
gives $H^4(A_{X_q};R)\cong H_Q^4(X_q;R)\cong C_{p^e}$. Theorem~\ref{james:all-powers} and
Corollary~\ref{james:entire-tower} now apply.
\end{proof}

\begin{corollary}[Unbounded nilpotence indices at fixed type]\label{symplectic:unbounded}
For each prime $p$, there is a family of finite connected quandles of fixed type for which the
nilpotence indices, under the cup product, of the degree-four IK images of the generators
in Theorem~\ref{binary:main} are unbounded over $\Z_{(p)}$. The type can be taken to be $3$ for $p=2$ and $2$ for odd $p$. Moreover, for any
fixed $m\ge1$, the groups $H_{4m-1}^Q(X_q;\Z_{(p)})$ in this family contain cyclic direct
summands of unbounded order.
\end{corollary}
\begin{proof}
For $j\ge0$, set
\[
 q_j\coloneqq \begin{cases}3^{2^j},&p=2,\\2^{(p-1)p^j},&p\text{ odd}.\end{cases}
\]
Fermat's little theorem and the binomial expansion give $e_j\coloneqq \nu_p(q_j^2-1)=j+3$ for
$p=2$ and $e_j=j+\nu_p(2^{2(p-1)}-1)$ for odd $p$. By \eqref{symplectic:power-order}, the
corresponding $y_R$ has nilpotence index $N_p(e_j)$, which is unbounded as $e_j\to\infty$. For
fixed $m$, we have $\nu_p(m!)<e_j$ for all sufficiently large $j$.
Corollary~\ref{james:high-torsion} then gives the direct summand $C_{p^{e_j}}$.
\end{proof}
By \eqref{james:coefficient-power-order}, for $s\ge1$ the reduction of $y_R$ modulo $p^s$ has
nilpotence index $N_p(\min(s,e_j))\le N_p(s)$. The degree-three quandle homology for odd
$q\notin\{3,5,7,9,27\}$ was computed in \cite[Theorem 8.6(I)]{NosakaHomotopical}. For $q>27$ the
existence of the cyclic direct summands also follows from the additive calculation in all degrees with the field characteristic inverted \cite[Example 9.2]{Model}. For the family $q=q_j$
above, Corollary~\ref{james:high-torsion} identifies the divided power $\gamma_m(y_R)$ as a
generator of a cyclic direct summand of order $p^{e_j}$ in $H_Q^{4m}(X_{q_j};R)$ whenever
$\nu_p(m!)<e_j$.

We now compare the operation values in Theorem~\ref{binary:main} with those of the ordinary
Steenrod operations.

\begin{proposition}\label{binary:squares}
In the setting of Theorem~\ref{binary:main}, if $e\ge2$, all positive-degree operations in
the Steenrod algebra act trivially on $H_R^*(X_q;\F_p)$ and $H_Q^*(X_q;\F_p)$. If $e=1$, then $p$ is odd, all
reduced powers $P^i$ for $i>0$ vanish, and the first Bockstein
$\beta_1:H_Q^3(X_q;\F_p)\to H_Q^4(X_q;\F_p)$ is nonzero.
\end{proposition}
\begin{proof}
Use $G,H,h$ from Proposition~\ref{binary:quandle} and put
$\mathcal{K}\coloneqq \Omega A_{X_q,R}$. Let $(BG)^\wedge_p$ denote the Bousfield--Kan
$p$-completion of $BG$ \cite[Chapter I]{BousfieldKan}. By Lemma~\ref{good:marked-pinch},
$A_{X_q,R}\simeq L_R$. Since $G$ is finite and $p\nmid |H|$, Proposition~3.3 of \cite{Model}
gives $L_R\simeq (BG)^\wedge_p$.

For odd $p$, a Sylow $p$-subgroup of $G$ is cyclic, and its normalizer acts on it through the
order-two group generated by inversion \cite[proofs of Lemmas 3.4--3.5]{Hutchinson}. For $p=2$,
a Sylow $2$-subgroup is generalized quaternion \cite[proof of Corollary 3.3]{Hutchinson}. In the case $p=2$,
the field characteristic $\ell$ is odd. The conjugates of $h$, all of order $\ell$, generate
$G$, so every homomorphism $G\to C_2$ is trivial. In both cases Benson's calculation \cite[Section 1.16, Corollary 4.7.5(iii)]{Benson} gives
generators $\tau_2$ and $\xi_3$ of the Pontryagin algebra $H_*(\mathcal{K};\F_p)$ in degrees $2$
and $3$, respectively. We use only that these classes generate, which also covers the
exceptional product when $p^e=3$.

First consider the dual homology operations. For odd $p$, the operation $P_*^i$ for $i>0$ lowers
degree by $2i(p-1)\ge4$, so it vanishes on both generators. By \eqref{binary:target-homology}
and \eqref{james:criterion-loop-low} with $N=4$, we have
$H_1(\mathcal{K};R)=H_3(\mathcal{K};R)=0$ and $H_2(\mathcal{K};R)=C_{p^e}$. Hence
$H_1(\mathcal{K};\F_p)=0$ and $\beta_{1,*}\tau_2=0$. For $p=2$, using also
$\operatorname{Sq}^i1=0$ for $i>0$, the only operation on a generator that can remain nonzero is
$\operatorname{Sq}_*^1\xi_3=\beta_{1,*}\xi_3$.

The universal coefficient theorem gives isomorphisms
$H_3(\mathcal{K};D)\cong\Tor_R(C_{p^e},D)\cong D[p^e]$, natural in the $R$-module $D$, where
$D[p^e]\coloneqq \ker(p^e:D\to D)$. For $e\ge2$, the kernel equals $D$ for both $D=\Z/p^2$ and
$D=\F_p$, so $H_3(\mathcal{K};\Z/p^2)\to H_3(\mathcal{K};\F_p)$ is surjective. The homology
coefficient sequence now gives $\beta_{1,*}\xi_3=0$.

The Cartan formula for the Pontryagin product extends the vanishing on the generators to
all of $H_*(\mathcal{K};\F_p)$. Thus the duals of all positive-degree ordinary Steenrod operations
vanish when $e\ge2$, and $P_*^i$ vanishes for $i>0$ at every odd prime under consideration.
Suspend to $V=\Sigma \mathcal{K}_+$ and apply the Cartan formula for external products to each
word-length summand of the James stable splitting \eqref{loop:stable-words}.
Proposition~\ref{james:space} and the injectivity of $\pi_X^*$ transfer the vanishing to rack
and quandle cohomology.

Finally suppose $e=1$. For $p=2$, the integer $q$ is odd and $8\mid q^2-1$, so the case $e=1$
occurs only at an odd prime. Theorem~\ref{good:degree3} with $n=3$ gives $H_2^Q(X_q;R)=0$ and
$H_3^Q(X_q;R)\cong C_p$. The universal coefficient isomorphisms
$H_Q^3(X_q;D)\cong\Hom_R(C_p,D)\cong D[p]$ are natural in $D$. For $D=\Z/p^2$ we have $D[p]=pD$,
so reduction $H_Q^3(X_q;\Z/p^2)\to H_Q^3(X_q;\F_p)$ is zero. Hence the connecting homomorphism
$\beta_1$ of $0\to\F_p\to\Z/p^2\to\F_p\to0$ is injective on $H_Q^3(X_q;\F_p)\cong\F_p$.
\end{proof}

For $e\ge2$, all positive-degree operations in the Steenrod algebra act trivially. By
\eqref{symplectic:tower}, $C_{p,r,Q}$ takes a nonzero value on $y_R\bmod p^r$ for every $r\ge1$.

Table~\ref{tab:symplectic-examples} lists three cases of Theorem~\ref{binary:main}, with
$R=\Z_{(p)}$, $y_R=\Phi_R^4(a)$ for a generator $a$, and $|t|=4$. The Bockstein $\beta_1$
is the map from degree $3$ to degree $4$ in quandle cohomology over $\F_p$. The orders in the
last column hold for every $r\ge1$, with order $1$ denoting the zero class.

\begin{table}[htbp]
\centering
\caption{Cup powers, Bocksteins, and Pontryagin--Thomas values for three symplectic quandles.}
\label{tab:symplectic-examples}
\smallskip
\begin{tabular}{@{}ccccc@{}}
\toprule
$(q,p,e)$ & $|X_q|$ & $R[y_R]$ & $\beta_1$ &
$\operatorname{ord}C_{p,r,Q}(y_R\bmod p^r)$\\
\midrule
$(3,2,3)$ & $8$ & $R[t]/(8t,4t^2,t^4)$ & $0$ & $2^{\min(r,2)}$\\
$(4,3,1)$ & $15$ & $R[t]/(3t,t^3)$ & $\ne0$ & $1$\\
$(8,3,2)$ & $63$ & $R[t]/(9t,3t^3,t^6)$ & $0$ & $3$\\
\bottomrule
\end{tabular}
\end{table}
The last two examples in Table~\ref{tab:symplectic-examples} both have type $2$. The
values $e=1$ and $e=2$ give different behavior for the first Bockstein and for
$C_{p,r,Q}$ on reductions of $y_R$.

In the eight-element example, $G=\mathrm{SL}_2(\F_3)$ is the binary tetrahedral group,
denoted by $P_{24}$ in \cite[\S4.1]{TomodaZvengrowski}. The ring calculation of
Tomoda and Zvengrowski~\cite[Theorem 4.4]{TomodaZvengrowski}, localized at $2$, together with
$q_G^*$ from Lemma~\ref{good:marked-pinch}, gives $H^*(A_{X_3};R)=R[a]/(8a)$ with
$R=\Z_{(2)}$ and $|a|=4$. By \eqref{eq:pontryagin-lift},
$C_{2,r,A_{X_3}}(a\bmod2^r)=a^2\bmod2^{r+1}$. Since
$a^2$ generates $H^8(A_{X_3};R)\cong C_8$, its reduction has order $2^{\min(r+1,3)}$ by the
coefficient sequence at the start of Section~\ref{james:section}. Its IK image is zero by
\eqref{symplectic:tower}. For $r=1$, the order of $a^2\bmod4$ is $4$, and the order of
$C_{2,1,Q}(y_R\bmod2)=y_R^2\bmod4$ is $2$.

For $(q,p)=(9,2)$, Theorem~\ref{binary:main} gives $H_Q^4(X_9;R)=\langle y_R\rangle\cong C_{16}$.
Thus $2y_R$ has order $8$ and is already nonzero modulo $4$. By \eqref{eq:pontryagin-lift}
and \eqref{james:coefficient-power-order} with $b=1$, $m=2$ and $s=r+1$, the class
$C_{2,r,Q}(2y_R\bmod 2^r)$ vanishes for $r\le2$ and has order $2$ for every $r\ge3$.

\subsection{Degree-six classes of cyclic Alexander quandles}

For an application of Theorem~\ref{james:all-powers} with $N=6$, put $R\coloneqq \Z_{(7)}$ and
for $e\ge1$ take
\[
 Q_e\coloneqq \Z/7^e,\qquad x*y\coloneqq T_ex+(1-T_e)y,\qquad
 T_e^2+T_e+1=0,\quad T_e\equiv2\pmod7\,.
\]
Hensel's lemma gives a unique $T_e$ with the stated residue. It is a unit of order $3$. The
additive quandle homology of
this family is covered by \cite[Theorem 1.1]{NosakaPrime} for $e=1$ and by \cite[Theorem
A]{Model} for general $e$.

\begin{example}\label{cyclic:sixth-degree}
The quandle $Q_e$ is finite and connected of type $3$, with
\[
 H_j(A_{Q_e};R)=0\quad(1\le j\le4),\qquad
 H_5(A_{Q_e};R)\cong C_{7^e}\,.
\]
Thus $\Phi_R^6:H^6(A_{Q_e};R)\to H_Q^6(Q_e;R)$ is an isomorphism between cyclic groups of order
$7^e$. Theorem~\ref{james:all-powers} and Corollary~\ref{james:entire-tower} apply to the image
of a generator with $N=6$. In particular, the cup-product nilpotence indices of the specified
degree-six IK images are unbounded within this family of type $3$.
\end{example}
\begin{proof}
Write $\tau_a(x)=x+a$ and $\sigma(x)=T_ex$. In the right action convention the right
translation at $y$ is $\eta(y)=\sigma\tau_{(1-T_e)y}$. Since $\eta(0)=\sigma$ and $1-T_e$
is a unit, these maps generate $G=C_{7^e}\rtimes\langle\sigma\rangle$. Thus $Q_e$ is
connected with stabilizer $H=\langle\sigma\rangle\cong C_3$ at $0$. Since
$\eta(y)=\tau_{-y}\sigma\tau_y$, every right translation has order $3$.

The periodic resolution of the cyclic group \cite[(I.6.3)]{Brown} gives
$H_{2k-1}(C_{7^e};R)=C_{7^e}$ and zero homology in positive even degrees. Since
$\sigma^{-1}\tau_a\sigma=\tau_{T_ea}$, conjugation by $\sigma^{-1}$ induces multiplication by
$T_e$ in degree one. The universal coefficient theorem and the $k$th power of the degree-two
cohomology generator show that the action in degree $2k-1$ is multiplication by $T_e^k$
\cite[Chapter V, \S3, Exercises 3--4]{Brown}. Since $3$ is invertible in $R$, only the
coinvariants remain in the Lyndon--Hochschild--Serre spectral sequence \cite[Theorem
VII.6.3]{Brown}. If $3\nmid k$, then $T_e^k-1$ is a unit. Thus in positive degrees $H_j(G;R)$ is
$C_{7^e}$ for $j\equiv5\pmod6$ and zero otherwise. Lemma~\ref{good:marked-pinch} and
Lemma~\ref{james:degree4} with $N=6$ give the asserted homology groups and the IK isomorphism.
The nilpotence index is $N_7(e)$ by Theorem~\ref{james:all-powers}, and it is unbounded as
$e\to\infty$.
\end{proof}
\appendix
\section{The relative chain model for the IK comparison}\label{appcomp:section}

We model $F$ on integral bar chains, keeping it constant on $BH$ and the cylinder. The prism
and cone paths give the signed map from the suspended fiber to the mapping cone in
Theorem~\ref{thm:comparison}.

Abbreviate $\mathcal N=H\backslash M$ of Section~\ref{sec:ik-setup} to $N$ and write the
integral cellular chain complexes as
\[
 P\coloneqq C_*^{\mathrm{cell}}(E),\qquad \mathcal A=C_*^{\mathrm{cell}}(M),\qquad
 S\coloneqq C_*^{\mathrm{cell}}(N),\qquad D\coloneqq C_*^{\mathrm{cell}}(A_X)\,.
\]
Here $q:P\to D$ is the orbit quotient, $\varepsilon$ is the augmentation, and $p_0\in P_0$ is
the chain of the vertex $x_0$. For suspension of chain complexes, we use
$(\susp C)_n\coloneqq C_{n-1}$ and $d(\susp c)=-\susp (dc)$.

For $x\in X$ and $z\in E$, define $T_x(t,z)$, $0\le t\le1$, by the path
$z\to x\to z^{\eta(x)}$. The first half
runs linearly in barycentric coordinates from $z$ to $x$ in the cone obtained by appending $x$
to the simplex. The second half reverses the corresponding path from $z^{\eta(x)}$ to $x$. The
identity $\eta(x^g)=g^{-1}\eta(x)g$ gives $T_x(t,z)^g=T_{x^g}(t,z^g)$. The right action of a point $[g;t_1,x_1,\ldots,t_n,x_n]\in M$ on $z$ starts at $z^g$
and applies $T_{x_i}(t_i,-)$ for $i=1,\ldots,n$ in that order. The endpoint identities
$T_{x_i}(0,-)=\id$ and $T_{x_i}(1,z)=z^{\eta(x_i)}$, together with equivariance, give
compatibility with the face identifications in \eqref{eq:rack-differential}. Equivariance
also allows the vertex of the second factor to pass to the left of the earlier labels,
as in the monoid multiplication. The resulting right action is cellular and makes $P$
a right differential graded $\mathcal A$-module.

Since $H$ fixes $x_0$, the map $m\mapsto x_0\cdot m$ factors through the covering quotient
$v:M\to N$ and defines an orbit map $e:N\to E$.
\begin{lemma}\label{appcomp:orbit-cancellation}
On cellular chains the orbit map induces the finite-difference formula \eqref{eq:Psi}.
\begin{equation}\label{appcomp:Psi}
 e_\#(r;\boldsymbol x)=\Psi_n(r;\boldsymbol x)\,.
\end{equation}
Moreover, for homogeneous chains $c\in P$ and $a\in\mathcal A$,
\begin{equation}\label{appcomp:chain-action-cancellation}
 q(c\cdot a)=\varepsilon(a)q(c)\,.
\end{equation}
In particular, $q\Psi_n=0$ for $n\ge1$.
\end{lemma}
\begin{proof}
Write $(c,x)$ for the simplex obtained by appending $x$ to $c$. Use the generating edge $E_x$ of
Section~\ref{sec:loop-structures}. With the product orientation on $c\times E_x$, the two cones
obtained by appending $x$ give
\[
 c\cdot E_x=(c,x)-(c^{\eta(x)},x)\,.
\]
Iterating this formula in coordinate order, starting from the vertex $r=x_0^g$, produces a
factor $1-R_i$ at each step and gives \eqref{appcomp:Psi}. Since $x^{\eta(x)}=x$, the two terms
in the formula for the generating edge lie in the same orbit of the diagonal $G$-action. Thus
$q(c\cdot E_x)=0$. On vertices, $q(c^g)=q(c)$, and every cube can be written as
$gE_{x_1}\cdots E_{x_n}$ by \eqref{loop:algebra-normal-form}. This proves
\eqref{appcomp:chain-action-cancellation}. Since $e_\#v_\#a=p_0\cdot a$ and every cube lifts to
the cover, we also get $q\Psi_{>0}=0$.
\end{proof}
For the IK map use the cone operator $h_0\coloneqq h_{x_0}$ that prepends $x_0$. By
Lemma~\ref{lem:finite-difference}, the map $\Theta:\susp \widetilde S\to\widetilde D$ given by
$\Theta(\susp c)=qh_0\Psi(c)$ is a chain map.

Use $Q\coloneqq B(E,M,*)$ to model $F$ on chains. The product cell structure on the bar
simplices
and their factors gives the normalized two-sided bar complex
$T\coloneqq B(P,\mathcal A,\Z)$, with degeneracies given by insertion of the unit. In the
internal differential on $p[a_1|\cdots|a_k]$, the tensor-product signs use the suspended
degree of each $a_i$. Since $E$ is contractible, the projection to $BM$ is an equivalence
in each simplicial degree. The degeneracy maps are inclusions of subcomplexes, so
$Q\to BM$ is an equivalence on realization. Define
\begin{equation}\label{appcomp:beta}
 \beta:T\longrightarrow D,\qquad
 \beta(p[])\coloneqq q(p),\qquad \beta(p[a_1|\cdots|a_k])\coloneqq 0\quad(k>0)\,.
\end{equation}
The two end faces of a bar of length one have equal images under $\beta$ and opposite signs by
\eqref{appcomp:chain-action-cancellation}, so $\beta$ is a chain map.

\begin{lemma}\label{appcomp:relative-beta}
Under the comparison preserving the marking and the prescribed maps on $BH$ and on the
cylinder, $\beta$ represents $F:Y\to A_X$ on chains.
\end{lemma}
\begin{proof}
Let $Y'$ be the double mapping cylinder obtained from $Y$ by replacing the $BG$ end by
$B(E,G,*)$ and the middle space of the span by
$B(H\times\mathbb N)$. The map to the left end is the composite of
$B\operatorname{pr}_1:B(H\times\mathbb N)\to BH$ and the map $BH\to B(E,G,*)$ given by the
$H$-fixed point $x_0$ and the inclusion $H\hookrightarrow G$. The map to the right end is
induced by
$(u,m)\mapsto uh^m$. Since $E$ is contractible and $B\mathbb N\simeq B\Z$ \cite[Lemma
3.2]{Mostovoy}, projection and inclusion give an equivalence $\nu:Y'\to Y$ that is the identity on $BH$.

The edge $w(t)=[1;t,x_0]:1\to h$ commutes with $H$ and fixes $x_0$. Thus
$f_t(u,m)\coloneqq uw(t)^m$ is a monoid homomorphism to $M_0\coloneqq B(H;\{x_0\})$. Define
$J:Y'\to Q$ by gluing the inclusion of $B(E,G,*)$, the maps $Bf_t$ along the cylinder,
and the map on $BH$ determined by $x_0$ and $H\hookrightarrow M$. In fact, $f_0(u,m)=u$ and
$f_1(u,m)=uh^m$. Since $f_t(u,0)=u$, the map $J$ is constant in the interval direction on the
cylinder $BH\times[0,1]$ over the section $u\mapsto(u,0)$. The images of $BH$ and the cylinder lie in $Q_0\coloneqq B(\{x_0\},M_0,*)$.

Write $\mathcal E$ for the realization of Mostovoy's simplicial group. After group completion
$M\to\mathcal E$, the same formula extends $f_t$ to $m\in\Z$. When the cylinder in the proof of
\cite[Theorem 2.1, (2.6)--(2.7)]{Model} is traversed from the $BG$ end to the $BH$ end, the
homomorphism on $H\times\Z$ is also the inclusion on $H$ and sends $(1,1)$ to the image of
$w(t):1\to h$. Thus $Y'\to Q\to BM\to B\mathcal E$ is the composite of $\nu$ with the marked
equivalence of that reference. Since $Q\to BM\to B\mathcal E$ is also an equivalence, $J$ is a
homotopy equivalence \cite[Theorem 3.1]{Mostovoy}.

Define $F':Y'\to A_X$ by $[z;g_1|\cdots|g_k]\mapsto q(z)$ on $B(E,G,*)$ and by the
constant value $q(x_0)$ on the cylinder and the right end. The first formula is well defined
since $q(z^g)=q(z)$, and its value on the attaching space is $q(x_0)$. On the
left end, $F'_\#=\beta J_\#$ since the map forgets the bar coordinates. Take a cellular
approximation of the cylinder map in $Q_0$ relative to its ends. Since $\beta$ vanishes on
positive-degree chains in $Q_0$, the comparison extends over $Y'$ relative to $BH$.

In the model $(EG\times E)/G$ for the left end, the $EG$ coordinate of an edge $[z;g]$ is
$(1,g^{-1})$, and right translation of its zeroth face by $g$ gives $z^g$. In these coordinates
$F'$ and $F\nu$ are induced by the functors sending a vertex $(g,x)$ to $x$ and $x_0^g$,
respectively. The unique arrow $x\to x_0^g$ in the pair groupoid gives a $G$-equivariant natural
transformation that is the identity on $EH\times\{x_0\}$ and its $G$-translates. Extend the
resulting homotopy constantly over the cylinder and $BH$. Together with the chain
comparison above this gives
\begin{equation}\label{appcomp:beta-J}
 \beta J_\#\simeq F'_\#\simeq(F\nu)_\#\qquad\text{relative to }BH\,.
\end{equation}
The map $J$ identifies the two copies of $BH$, whose inclusions are cofibrations. A homotopy
inverse of $J$ under $BH$ therefore transports the chain comparison to $Q$.
\end{proof}

To pass from this relative comparison to the mapping cone, write
$W\coloneqq B(\Z,\Z H,\mathcal A)$ and $U\coloneqq B(\Z,\Z H,\Z)$ for the cellular chains of
$B(*,H,M)$ and $BH$, respectively. For $b=[u_1|\cdots|u_k]$ of length $k$, with $u_i\in H$, the
internal differential of $b\otimes a$ is $(-1)^k b\otimes da$. Projection to the base induces
$\pi:W\to U$, $b\otimes a\mapsto\varepsilon(a)b$, and the map from $BH$ induces $j_H:U\to T$,
$b\mapsto p_0b$. The map $B(*,H,M)\to H\backslash M$ to the covering quotient induces
$\rho_H:W\to S$, given by $[]\otimes a\mapsto v_\#a$ in bar length zero and by zero in positive
bar lengths. Since $\mathcal A$ is a free left $\Z H$-module in each degree, the augmentation of
the bar resolution $\rho_H:W\to\Z\otimes_{\Z H}\mathcal A\cong S$ is an integral
quasi-isomorphism. Set $f\coloneqq \Psi\rho_H:W\to P$, and let $i:P\to T$ be the inclusion in
bar length zero.

Define the degree-one prism operator $\mathsf P:W\to T$ by
\begin{equation}\label{appcomp:prism}
 \mathsf P([u_1|\cdots|u_k]\otimes a)\coloneqq (-1)^k p_0[u_1|\cdots|u_k|a]\,.
\end{equation}

\begin{lemma}\label{appcomp:prism-identity}
The operator $\mathsf P$ is the prism chain homotopy from $j_H\pi$ to $if$. It satisfies
\begin{equation}\label{appcomp:prism-homotopy}
 d_T\mathsf P+\mathsf P d_W=if-j_H\pi,\qquad \beta\mathsf P=0\,.
\end{equation}
\end{lemma}
\begin{proof}
The functor inducing $j_H\pi$ sends the arrow $u:um\to m$ of the action category, for $m\in M$,
to $u:x_0\to x_0$. The functor inducing $if$ sends it to the identity of $x_0m$. Since
$x_0u=x_0$, the arrow $m:x_0\to x_0m$ gives a natural transformation from the first functor to
the second. In the prism sum on the nerve, every term except the last insertion contains an
identity arrow and is degenerate. The normalized prism is therefore \eqref{appcomp:prism}, with
sign $(-1)^k$ for the last insertion. The boundary formula gives the terminal map minus the initial map, $if-j_H\pi$. For $k=0$ this is $i(p_0\cdot a-\varepsilon(a)p_0)$. Since $\mathsf P$
raises bar length by one, the two terms involving $da$ have opposite signs and cancel. The image
of $\mathsf P$ has positive bar length, so $\beta\mathsf P=0$.
\end{proof}

Pass to reduced chains at the basepoints and put $C\coloneqq ih_0f$ and
$L\coloneqq C-\mathsf P$. The identities $dh_0+h_0d=\id-p_0\varepsilon$ and
$\varepsilon f|_{\widetilde W}=0$ give $d_TC+Cd_W=if$. Hence \eqref{appcomp:prism-homotopy}
gives
\begin{equation*}
 d_TL+Ld_W=j_H\pi\,.
\end{equation*}
Equip $\operatorname{Cone}(j_H)\coloneqq \widetilde T\oplus \susp \widetilde U$ with
differential $d(t,\susp u)\coloneqq (d_Tt+j_Hu,-\susp d_Uu)$. Then the map
\begin{equation}\label{appcomp:I}
 I:\susp \widetilde W\longrightarrow\operatorname{Cone}(j_H),\qquad I(\susp w)\coloneqq (Lw,-\susp \pi w)
\end{equation}
is a chain map, since $dI(\susp w)=(-Ldw,\susp \pi dw)=I(-\susp dw)$.

For each point of $B(*,H,M)$, concatenate the cone path from the basepoint to its image in
$E$ with the reverse of the prism path from $BH$.
For $z=[(t_0,\ldots,t_n);z_0,\ldots,z_n]$ the cone path is
\[
 \tau\longmapsto [(1-\tau,\tau t_0,\ldots,\tau t_n);x_0,z_0,\ldots,z_n],
 \qquad 0\le\tau\le1.
\]
The concatenated paths define a map $B(*,H,M)\to\hofib(BH\to Q)$ over $BH$.
With the interval oriented first, the cone path induces $h_0$ on chains. The reversed
prism contributes $-\mathsf P$, so the concatenation induces $L=ih_0f-\mathsf P$. For a point $m\in M$ in the
fiber, the projection $Q\to BM$ sends the cone path to a constant path and the prism to the bar
edge $\ell_m$. The map on the fiber is therefore the bar--loop equivalence of the connected
monoid $M$ followed by loop reversal \cite[Proposition 4.17(1), Theorem 5.6]{Vogt}. Comparison
of the fibrations over $BH$ gives an equivalence on their total spaces. Via $B(*,H,M)\simeq N$
and the marked comparison $c_M:M\to\Omega Y$, we obtain
\[
 \kappa v\simeq\operatorname{inv}_{\Omega A_X}\circ(\Omega F)\circ c_M\,.
\]
Appending the path from the $BH$ end to the cone point gives
$\Sigma\hofib(BH\to Q)\to\Cofib(BH\to Q)$. With the mapping-cone differential above, $\susp u$
denotes the cone cell oriented from the cone point toward the base cell $u$. The return path has
the opposite orientation, so its chain component is $-\susp \pi w$, as in \eqref{appcomp:I}.

\section*{Acknowledgments}
The author used OpenAI's ChatGPT to write code for checking the finite-quandle examples and discussions in Appendix A, and to revise the exposition of the proofs. The author takes responsibility for the content of the paper, including the
interpretation of the computations.

\begingroup
\linespread{1}\selectfont

\endgroup
\smallskip
\noindent Department of Mathematics, School of Science, Institute of Science Tokyo

\begin{thebibliography}{CFLM}
\small
\raggedright
\setlength{\itemsep}{0.4pt}
\setlength{\parsep}{0pt}
\setlength{\parskip}{0pt}
\interlinepenalty=10000

\bibitem[AC]{AdamsonComparison}
J.~A. Arciniega-Nev\'arez, J.~L. Cisneros-Molina and L.~J. S\'anchez Salda\~na,
\emph{Relative group homology theories with coefficients and the comparison homomorphism},
Quaest. Math. \textbf{44} (2021), 1107--1132.
\href{https://doi.org/10.2989/16073606.2020.1770895}{doi:10.2989/16073606.2020.1770895};
\href{https://arxiv.org/abs/1809.01209v2}{arXiv:1809.01209v2}.

\bibitem[Be]{Benson}
D.~J. Benson,
\emph{Classifying spaces of finite groups of tame representation type},
\href{https://arxiv.org/abs/2208.07913v2}{arXiv:2208.07913v2}, 2025.

\bibitem[BK]{BousfieldKan}
A.~K. Bousfield and D.~M. Kan,
\emph{Homotopy Limits, Completions and Localizations},
Lecture Notes in Math. 304, Springer-Verlag, Berlin--New York, 1972.
\href{https://doi.org/10.1007/978-3-540-38117-4}{doi:10.1007/978-3-540-38117-4}.

\bibitem[BT]{BrowderThomas}
W. Browder and E. Thomas,
\emph{Axioms for the generalized Pontryagin cohomology operations},
Quart. J. Math. Oxford (2) \textbf{13} (1962), 55--60.
\href{https://doi.org/10.1093/qmath/13.1.55}{doi:10.1093/qmath/13.1.55}.

\bibitem[Br]{Brown}
K.~S. Brown,
\emph{Cohomology of Groups},
Grad. Texts in Math. 87, Springer-Verlag, New York, 1982.
\href{https://doi.org/10.1007/978-1-4684-9327-6}{doi:10.1007/978-1-4684-9327-6}.

\bibitem[Cl]{Clauwens}
F.~J.-B.~J. Clauwens,
\emph{The algebra of rack and quandle cohomology},
J. Knot Theory Ramifications \textbf{20} (2011), no.~11, 1487--1535.
\href{https://doi.org/10.1142/S0218216511010073}{doi:10.1142/S0218216511010073};
\href{https://arxiv.org/abs/1004.4423}{arXiv:1004.4423}.

\bibitem[CFLM]{CFLM}
S. Covez, M.~A. Farinati, V. Lebed and D. Manchon,
\emph{Bialgebraic approach to rack cohomology},
Algebr. Geom. Topol. \textbf{23} (2023), no.~4, 1551--1582.
\href{https://doi.org/10.2140/agt.2023.23.1551}{doi:10.2140/agt.2023.23.1551};
\href{https://arxiv.org/abs/1905.02754}{arXiv:1905.02754}.

\bibitem[DH]{DevalapurkarHaine}
S. Devalapurkar and P. Haine,
\emph{On the James and Hilton--Milnor splittings, and the metastable EHP sequence},
Doc. Math. \textbf{26} (2021), 1423--1464.
\href{https://doi.org/10.4171/DM/845}{doi:10.4171/DM/845}.

\bibitem[EG]{EtingofGrana}
P. Etingof and M. Gra\~na,
\emph{On rack cohomology},
J. Pure Appl. Algebra \textbf{177} (2003), no.~1, 49--59.
\href{https://arxiv.org/abs/math/0201290}{arXiv:math/0201290}.

\bibitem[Do]{Dokas}
I. Dokas,
\emph{Zinbiel algebras and commutative algebras with divided powers},
Glasgow Math. J. \textbf{52} (2010), 303--313.
\href{https://doi.org/10.1017/S0017089509990358}{doi:10.1017/S0017089509990358}.

\bibitem[FRS]{FRS}
R. Fenn, C. Rourke and B. Sanderson,
\emph{The rack space},
Trans. Amer. Math. Soc. \textbf{359} (2007), no.~2, 701--740.
\href{https://doi.org/10.1090/S0002-9947-06-03912-2}{doi:10.1090/S0002-9947-06-03912-2};
\href{https://arxiv.org/abs/math/0304228v2}{arXiv:math/0304228v2}.

\bibitem[GJ]{GoerssJardine}
P.~G. Goerss and J.~F. Jardine,
\emph{Simplicial Homotopy Theory},
Modern Birkh\"auser Classics, Birkh\"auser, Basel, 2009.
\href{https://doi.org/10.1007/978-3-0346-0189-4}{doi:10.1007/978-3-0346-0189-4}.

\bibitem[Hu]{Hutchinson}
K. Hutchinson,
\emph{A Bloch--Wigner complex for $\mathrm{SL}_2$},
J. K-Theory \textbf{12} (2013), no.~1, 15--68.
\href{https://doi.org/10.1017/IS013003031JKT13222}{doi:10.1017/IS013003031JKT13222};
\href{https://arxiv.org/abs/1107.0264v2}{arXiv:1107.0264v2}.

\bibitem[In]{Inoue2026}
A. Inoue,
\emph{Quandle homology and relative group homology},
\href{https://arxiv.org/abs/2607.00701v2}{arXiv:2607.00701v2}, 2026.

\bibitem[IK]{IK}
A. Inoue and Y. Kabaya,
\emph{Quandle homology and complex volume},
Geom. Dedicata \textbf{171} (2014), 265--292.
\href{https://doi.org/10.1007/s10711-013-9898-2}{doi:10.1007/s10711-013-9898-2};
\href{https://arxiv.org/abs/1012.2923v3}{arXiv:1012.2923v3}.

\bibitem[IT]{IT}
K. Ishikawa and K. Tanaka,
\emph{Quandle colorings vs. biquandle colorings},
Topology Appl. \textbf{345} (2024), Paper No.~108832.
\href{https://doi.org/10.1016/j.topol.2024.108832}{doi:10.1016/j.topol.2024.108832};
\href{https://arxiv.org/abs/1912.12917v2}{arXiv:1912.12917v2}.

\bibitem[Ja]{James}
I.~M. James,
\emph{Reduced product spaces},
Ann. of Math. (2) \textbf{62} (1955), no.~1, 170--197.
\href{https://doi.org/10.2307/2007107}{doi:10.2307/2007107};
\href{https://webhomes.maths.ed.ac.uk/~v1ranick/papers/jamesred.pdf}{PDF of the original article}.

\bibitem[LN]{LN}
R.~A. Litherland and S. Nelson,
\emph{The Betti numbers of some finite racks},
J. Pure Appl. Algebra \textbf{178} (2003), no.~2, 187--202.
\href{https://doi.org/10.1016/S0022-4049(02)00211-6}{doi:10.1016/S0022-4049(02)00211-6};
\href{https://arxiv.org/abs/math/0106165}{arXiv:math/0106165}.

\bibitem[Lo]{Loday}
J.-L. Loday,
\emph{Cup-product for Leibniz cohomology and dual Leibniz algebras},
Math. Scand. \textbf{77} (1995), no.~2, 189--196.
\href{https://doi.org/10.7146/math.scand.a-12560}{doi:10.7146/math.scand.a-12560}.

\bibitem[Mi]{Milnor}
J.~W. Milnor,
\emph{On the construction $FK$},
in J.~F. Adams, \emph{Algebraic Topology: A Student's Guide},
London Math. Soc. Lecture Note Ser. 4, Cambridge Univ. Press, 1972, 118--136.
\href{https://doi.org/10.1017/CBO9780511662584.011}{doi:10.1017/CBO9780511662584.011}.

\bibitem[Mo]{Mostovoy}
J. Mostovoy,
\emph{Pre-crossed modules and rack homology},
Commun. Math. \textbf{30} (2022), no.~2, 119--128.
\href{https://doi.org/10.46298/cm.10153}{doi:10.46298/cm.10153};
\href{https://arxiv.org/abs/2210.02708v2}{arXiv:2210.02708v2}.

\bibitem[NN]{NavasNelson}
E.~A. Navas and S. Nelson,
\emph{On symplectic quandles},
Osaka J. Math. \textbf{45} (2008), no.~4, 973--985.
\href{https://arxiv.org/abs/math/0703727v2}{arXiv:math/0703727v2}.

\bibitem[NoH]{NosakaHomotopical}
T. Nosaka,
\emph{Homotopical interpretation of link invariants from finite quandles},
Topology Appl. \textbf{193} (2015), 1--30.
\href{https://doi.org/10.1016/j.topol.2015.05.087}{doi:10.1016/j.topol.2015.05.087};
\href{https://arxiv.org/abs/1210.6528v6}{arXiv:1210.6528v6}.

\bibitem[NoP]{NosakaPrime}
T. Nosaka,
\emph{On quandle homology groups of Alexander quandles of prime order},
Trans. Amer. Math. Soc. \textbf{365} (2013), no.~7, 3413--3436.
\href{https://doi.org/10.1090/S0002-9947-2013-05754-6}{doi:10.1090/S0002-9947-2013-05754-6}.

\bibitem[No2]{Model}
T. Nosaka,
\emph{Loop Space Models and Integral Homology of Connected Quandles},
\href{https://arxiv.org/abs/2610.38066v1}{arXiv:2610.38066v1}, 2026.

\bibitem[Se]{SerreCubical}
J.-P. Serre,
\emph{Homologie singuli\`ere des espaces fibr\'es. Applications},
Ann. of Math. (2) \textbf{54} (1951), no.~3, 425--505.
\href{https://doi.org/10.2307/1969485}{doi:10.2307/1969485}.

\bibitem[TZ]{TomodaZvengrowski}
S. Tomoda and P. Zvengrowski,
\emph{Remarks on the cohomology of finite fundamental groups of $3$-manifolds},
Geom. Topol. Monogr. \textbf{14} (2008), 519--556.
\href{https://doi.org/10.2140/gtm.2008.14.519}{doi:10.2140/gtm.2008.14.519};
\href{https://arxiv.org/abs/0904.1876}{arXiv:0904.1876}.

\bibitem[Vo]{Vogt}
R.~M. Vogt,
\emph{Homotopy homomorphisms and the classifying space functor},
J. Homotopy Relat. Struct. \textbf{10} (2015), 875--921.
\href{https://doi.org/10.1007/s40062-014-0087-2}{doi:10.1007/s40062-014-0087-2};
\href{https://arxiv.org/abs/1203.4978v4}{arXiv:1203.4978v4}.

\bibitem[Wu]{WuSimplicial}
J. Wu,
\emph{Simplicial Objects and Homotopy Groups},
in \emph{Braids}, Lect. Notes Ser. Inst. Math. Sci. Natl. Univ. Singap. 19,
World Scientific, Hackensack, NJ, 2010, 31--181.
\href{https://doi.org/10.1142/9789814291415_0002}{doi:10.1142/9789814291415\_0002};
\href{https://imsarchives.nus.edu.sg/oldwww/Programs/braids/files/JieWuSimplicial.pdf}{Open-access PDF}.

\end{thebibliography}
\end{document}